\documentclass[11pt,reqno]{article}

\usepackage{amssymb,amsmath,amsthm,mathtools,wasysym,calc,verbatim,enumitem,tikz,pgfplots,hyperref,url,mathrsfs,fullpage,bbm,comment}
\usepackage[noadjust]{cite}
\usepackage{authblk}
\usepackage{comment,fullpage}
\usepackage[normalem]{ulem}
\mathtoolsset{showonlyrefs}
\pgfplotsset{compat=1.18}
\usepackage{mathabx}

\newtheorem{theorem}{Theorem}

\newtheorem{proposition}[theorem]{Proposition}
\newtheorem{lemma}[theorem]{Lemma}
\newtheorem{claim}[theorem]{Claim}

\newtheorem{conj}[theorem]{Conjecture}
\newtheorem{fact}[theorem]{Fact}

\newtheorem*{claim*}{Claim}

\theoremstyle{remark}
\newtheorem*{remark*}{Remark}

\usepackage{thmtools, thm-restate}

\numberwithin{theorem}{section}

\DeclarePairedDelimiter\norm{\lVert}{\rVert}

\renewcommand{\phi}{\varphi}

\renewcommand{\leq}{\le}
\renewcommand{\geq}{\ge}

\newcommand{\eps}{\varepsilon}

\newcommand{\cL}{\mathcal L}

\newcommand{\cH}{\mathcal H}

\usepackage{cancel}

\newcommand{\cS}{\mathcal S}
\newcommand{\cP}{\mathcal P}

\newcommand{\E}{\mathbb{E}}

\def\1{\mathbbm{1}}

\newcommand{\Unif}{\operatorname{Unif}}

\usepackage{tikz-cd}
\usetikzlibrary{positioning,arrows.meta,calc}

\newcommand{\sol}[1]{%
  \mkern1mu\overline{\mkern-1mu #1 \mkern-1mu}\mkern1mu}

\renewcommand{\le}{\leqslant}
\renewcommand{\ge}{\geqslant}
\renewcommand{\P}{\mathbb{P}}

\newcommand{\R}{\mathbb R}

\newcommand{\bX}{\mathbf X}
\newcommand{\bY}{\mathbf Y}

\newcommand{\PP}{\mathbb P}

\newcommand{\bx}{\mathbf x}
\newcommand{\by}{\mathbf y}

\newcommand{\x}{\times}

\usepackage{thm-restate}
\newcommand{\Norm}[1]{\left\lVert#1\right\rVert}

\usepackage{array}



\DeclareMathOperator*{\argmin}{arg\,min}
\usepackage{esvect}

\title{Local-to-Global Theorems Based on Cones \\with Applications to Graph Coloring}
\author[1]{Anqi Li}
\author[1]{Shayan Oveis Gharan}
 
\affil[1]{University of Washington, \textsf{\{aqli,shayan\}@cs.washington.edu}}
\date{}

\begin{document}
\maketitle
\begin{abstract}
Local-to-Global theorems are a cornerstone of many of recent developments in theoretical computer science from coding theory and PCPs all the way to analyzing mixing time of Markov chains. Building on trickledown theorems on simplicial complexes, several techniques have been established to obtain tighter local-to-global theorems to analyze the spectral gap of the Glauber dynamics, namely matrix trickledown, ${\cal C}$-Lorentzian polynomials and more  recently Bochner type inequalities. 

In this paper, we adopt non-local local-to-global/Bochner type inequalities on ``cones'', corresponding to the local neighborhood of a given vertex, to show that the Glauber dynamics mixes in polynomial time to sample proper colorings of maximum degree $\Delta$ graphs as long as $q\geq 1.7408\Delta$.
\end{abstract}

\section{Introduction}
Let $G=(V,E)$ be a graph with $n$ vertices and maximum degree $\Delta>1$ and a family of list of colors $\{\cL_i\}_{i\in G}$. A proper coloring is an assignment $\bx\in \Pi_i \cL_i$ such that for any two adjacent vertices $i\sim j$, we have $x_i\neq x_j$. Let $\mu$ be the uniform distribution over all proper colorings. One can run a simple Markov chain to generate a sample from $\mu$: Start with an arbitrary proper coloring $\bx$, each time choose a uniformly random vertex $i$ and re-sample its color given the colors of its neighbors. 
Let $\frac{1}{n}\cP$ be the transition probability matrix of the Glauber dynamics Markov chain (the factor $1/n$ is for notational convenience as will be clear later; essentially it implies every vertex is chosen at a ``unit rate") then the mixing time is defined to be
$$ \max_{\bx} \min \bigg\{t: \bigg\|1_\bx^\top \bigg(\frac1{n}\cP \bigg)^t  - \mu \bigg\|_{\mathrm TV}\leq \tfrac{1}{4} \bigg\}.$$
To bound the mixing time of $\frac1{n}\cP$ a method of choice is to bound its spectral gap. 
\begin{theorem}\label{thm:spec-gap}
If the second largest eigenvalue of $\frac{1}{n}P$ is at most $1-\eps$, or equivalently the second smallest eigenvalue of the Laplacian matrix $\frac{1}{n}D=I-\frac{1}{n}P$ is at least $\eps$, then the mixing time of $\frac{1}{n}\cP$ is at most $O\big(n\cdot \log(\frac{1}{\mu_{\min}})/\eps \big)$
where $\mu_{\min}=\min_x\mu(x)$.
\end{theorem}

Suppose that the size of every list is at least $q$. It is a long standing open problem that if $q\geq \Delta+2$, then the Glauber dynamics mentioned above mixes rapidly. 
\begin{conj}\label{conj:coloring}
    For any graph $G$ with maximum degree $\Delta$, if $q\geq \Delta+2$, then the Glauber dynamics mixes in time polynomial in $n$. 
\end{conj}
In the early 1990s, Jerrum \cite{J95} showed that $q\geq 2\Delta+1$ suffices via a simple coupling approach. Vigoda \cite{V00} improved this to $11\Delta/6$ using a more sophisticated coupling. Later in a sequence of works \cite{CDMPP19, JMV26, CV25} the bound was improved to $(11/6 - \eps)\Delta$ using various more involved coupling techniques. 

\paragraph{Local to Global Theorems.} With the introduction of the spectral independence technique \cite{ALO24} a new challenge started: To come up with new local-to-global arguments to analyze Markov chains. It is worth noting that even though classical proof techniques such as the correlation decay or path coupling methods could be translated to spectral independence \cite{ALO24,CLV23,L21,BCCPSV22} it seemed the main strength of the spectral independence machinery stems from the local-to-global technique known as the trickledown theorem \cite{Opp18} which was used to count bases of  matroids \cite{ALOV24}. 

Towards this goal, analyzing Glauber dynamics for graph coloring was considered as a benchmark.  The first improvement over the classical trickledown technique came from the Matrix trickledown theorems \cite{ALO21,AKV24,CCCYZ25} which were used to obtain near optimal bounds for edge coloring problems \cite{WZZ24}, or to sample from the Sherrington-Kirkpatrick model \cite{AKV24}. Later, this was extended to the ${\cal C}$-Lorentzian polynomial framework which led to a ``spectral" Dobrushin's criterion \cite{LLO25, LO25} that itself lead to improved analysis of the multi-state hardcore model or improved analysis of the Hard core model on planar graphs \cite{GZ26}. Later on, G\"{o}bel, Jenssen, Michelen, Pappik, Perkins and Schiller \cite{GJMPPS26} observed that the matrix trickledown of \cite{CCCYZ25} can be seen as a discrete Bochner inequality. This led to several improvement for various problems \cite{CL26a,BLO26,CL26b}. In this paper, we establish the discrete Bochner inequality over $O(\Delta)$-sized subgraphs and then use it to analyze the Glauber dynamics for sampling proper colorings of graphs.

\paragraph{Background on Bochner Inequality.} The discrete Bochner type inequalities were recently  recently used in \cite{GJMPPS26} to re-prove the trickledown theorem of Chen, Chen, Chen, Yin, and Zhang \cite{CCCYZ25}. For background on this technique see \cite[Equation 1.16.6, Proposition 4.8.3]{BGL14} and references therein. Bochner-type methods have also been used to establish spectral-gap and
mixing-time bounds in a variety of settings, including the
Sherrington-Kirkpatrick model \cite{Wan26,BLO26}, matroid intersection \cite{CL26b}, and several recent works by OpenAI. 

Let $P_if$ be the function which averages $f$ on $\cL_i$ given the color of all other vertices, namely:
$$ P_if(\bx)=\E_\mu[f(X_1,\dots,X_n) | X_j=x_j,\forall j\neq i].$$
Throughout the paper we use bold letters to denote vectors. 

Then, we can write $P=\sum_i P_i$. In other words the $P_i$ operators give a decomposition of $P$ into projection matrices; and in most recent applications this has been perhaps the main starting point of Bochner type inequalities. 
Then, it follows that the transition probability matrix of the Glauber dynamics is $\frac1n{P}.$ We let $D_i=I-P_i$ and $D =\sum_i D_i.$
It is not hard to see that the operators $P_i,D_i$ are projections, i.e., for any function $f$, $P_i^2f=P_if$. Define the inner-product 
$$ \langle f,g\rangle_\mu=\E_{X\sim\mu}[f(X)g(X)].$$
It is not hard to see that $P_i,D_i$ are self-adjoint with respect to this inner-product, i.e., for two functions $f,g$, $\langle f, P_i g \rangle = \langle P_if , g \rangle$,
and since $D_i = I-P_i$, the self-adjointness of $D_i$ follows.
\begin{fact}\label{fact:LLvsL}
    If for any function $f:\Omega \to \R$,
    \begin{equation}\label{eq:DDvsDeps} \langle D f,D f\rangle_\mu\geq \eps\cdot\langle D f,f\rangle_\mu=\eps\cdot\sum_i \langle D_if,f\rangle_\mu,
    \end{equation}
    then the spectral gap of the Glauber dynamics $\frac1{n}\cP$ is at least $\eps/n.$
\end{fact}
\begin{proof}
Given an arbitrary function $g$, since $D\succeq 0$, we let $f=D^{-1/2} g$ (to be precise, we need to take the pseudo-inverse of $D$ since $D$ has a 0 eigenvalue). 
Having this,
    $$ \langle D g,g\rangle=\langle D f,D f\rangle\geq\eps\cdot\langle D f,f\rangle = \eps \cdot \langle g,g\rangle,$$
    where all the inner-products are with respect to $\mu$. This concludes the proof.
\end{proof}

\paragraph{Top Link Spectral Expansion.}
A natural approach to prove \eqref{eq:DDvsDeps} is to lower-bound $\langle D_if,D_jf\rangle_\mu$ for $i,j\in G$. This is very analogous to analyzing links of co-dimension 2 in simplicial complexes, see e.g., \cite{ALOV24}. See the following lemma and the proof. 
\begin{lemma}\label{lem:mixed-form-lower-bound}
Given a graph $H$ with vertices $1$ and $2$. If $H$ has no edges then $\langle D_1f,D_2f\rangle_\mu \geq 0$. Otherwise, suppose $s \geq 5$ and  $|\cL_1|, |\cL_2| \geq s+1$. Then for every $f \in L^2(\mu)$, we have 
\begin{equation}\label{eq:mixed-form-lower-bound}
    2 \langle D_1f,D_2f\rangle_\mu
    \geq
    -\frac{1}{s}\big(\norm{D_1f}_\mu^2+\norm{D_2f}_\mu^2\big).
\end{equation}
\end{lemma}
\begin{proof}
If $H$ has no edges then $\mu$ is a product measure, so
$D_1,D_2$ are commuting orthogonal projections: 
\[
    \langle D_1f,D_2f\rangle_\mu
    =\langle f,D_1D_2f\rangle_\mu
    =\norm{D_1D_2f}_\mu^2\geq0.
\]

Now suppose $H$ has an edge. We let $P$ be the simple
random walk on the bipartite graph with classes
$\{1\}\times\cL_1$ and $\{2\}\times\cL_2$, where
$(1,a)$ and $(2,b)$ are adjacent exactly when $a\ne b$. By \cite[Eq.~(4)]{CL26a},
\begin{equation}\label{eq:local-walk-bound}
    2\langle D_1f,D_2f\rangle_\mu
    \geq
    -\lambda_2(P)
    \bigl(\norm{D_1f}_\mu^2+\norm{D_2f}_\mu^2\bigr),
\end{equation}
where $\lambda_2(P)$ denotes the second largest eigenvalue. It follows from \cite[proof of Proposition~IV.4]{ALO21} that $\lambda_2(P)\leq \frac{1}{\sqrt{(\ell_1-1)(\ell_2-1)}}\leq 1/s$.

\end{proof}
Having this, we can obtain that the Glauber dynamics has a spectral gap at least $1/\Delta n$ if $q>2\Delta+2.$ In particular
\begin{align*}
    \langle Df,Df\rangle\geq \sum_i \langle D_if,f\rangle_{\mu} +\sum_{\{ij\}\in E(G)} 2\langle D_if,D_jf\rangle_\mu
    \geq \sum_{i} \bigg(1-\frac{\deg(i)}{q-\Delta} \bigg) \cdot \langle D_if,f\rangle_\mu \geq \frac{1}{\Delta}\langle Df,f\rangle_\mu.
\end{align*}

\subsection{Main Contribution}
The first idea to improve over the $q>2\Delta$ bound is to do a more careful analysis of $\langle D_if,D_jf\rangle_\mu$ for adjacent $i,j$. Here, it is insightful to go back to the simplicial complex view and analyze the link where every  vertex (other than $i,j$) is pinned to a color. The 1-skeleton of such a link corresponds to a complete bipartite graph where one  matching is removed; one side corresponds to all colors available to $i$ and the other side corresponds to all colors available to $j$. As we have seen in the proof of Lemma~\ref{lem:mixed-form-lower-bound}, the second eigenvalue of the simple random walk on this graph is $\leq 1/s$. But it turns out that the eigenvectors corresponding to this bad eigenvalue of $1/s$ assign opposite signs to vertices on the two sides of the bipartite graph. Now, imagine an odd cycle, say a triangle; it follows that the bad eigenvectors corresponding to the top links for the edges of the triangle are not consistent! So, a more careful analysis takes this info into account and exploits odd cycles. 
In particular, we prove a stronger version of the above lemma (analogous to \cite[Proposition IV.4]{ALO21}) 
$$ \langle \sol{D}_if,\sol{D}_jf\rangle \geq  \E_{\bx\sim\mu}\sum_{c\in \cL_i\cap \cL_j}(f(\bx)-f(\bx^{i\to c}))(f(\bx)-f(\bx^{j\to c}).$$
In the above $\sol{D}_i$ can be seen as an unnormalized Laplacian and $\bx^{i\to c}$ sets the color of vertex $i$ to $c$.
To put it differently,  say we abstract the info of such a link by defining a vector $\nabla_i f$ with coordinates $\nabla_i f(\bx)_c = f(\bx)-f(\bx^{i-\to c}),$ and similarly for $j$.  Then, instead of lower-bounding the cross term by a single scalar $1/s$, we lower-bound it by $\langle \sol{D}_i f,\sol{D}_j f\rangle \geq \langle \nabla_if,\nabla_jf\rangle$.
Having this, we use a matrix Hoeffding inequality to obtain the following.
\begin{restatable}{theorem}{mainlambdamin}\label{cor:main-spectral}
    If for every vertex $i$, $|\cL_i|\geq s+\deg(i)$ where $s > 2 |\lambda_{\min}(A_G)|$, then  the Glauber dynamics for proper list colorings of a graph with maximum degree $\Delta$ mixes in time $O(n^2)$.
\end{restatable}
Throughout the article, instead of lower-bounding $q-\Delta$, we lower-bound $s$; this allows our statements to extend naturally when we induct. As a special case if $G$ is the line graph of another graph, $\lambda_{\min}(A_G) \geq -2$ so the above bound shows that $s > 4$ is enough. 

The above theorem can be seen as analogous to the work of \cite{CCCYZ25} who showed that the Glauber dynamics mixes rapidly for the hardcore model as long as the activity threshold is $<\frac{1}{|\lambda_{\min(A_G)}|-1}$.
Interestingly, such a Markov chain may have an exponential mixing time above the tree uniqueness threshold, i.e., $\approx \frac{e}{\Delta}$ on $\Delta$-regular bipartite graphs (even if they are random \cite{DFJ99}). 
This may seem very surprising and possibly contradictory with Conjecture~\ref{conj:coloring}. 
The main observation is that if $G$ is (close to) bipartite the above argument leading to Theorem \ref{cor:main-spectral} fails as bad eigenvectors of links of co-dimension 2 ``match up" (as the graph has no odd cycles). So, the main way to make progress is to avoid doing a worst case analysis for links of co-dimension 2 which  only assumes every vertex has $s$-colors left (note that there is no analogue of $s$ for the hardcore model). 

In other words, the observation is that in graphs close to bipartite one typically expects many more than $s$ colors available to vertices under a random proper coloring. So, in our analysis we naturally have to analyze links of co-dimension more than 2 (in fact for similar reasons we need to analyze links of co-dimension  $\gtrsim \Delta$).  This was in a sense apparent to experts, but unfortunately, analyzing the matrix trickledown or ${\cal C}$-Lorentzian polynomials machinery on such links was difficult.

Instead, we analyzed such links using the new Bochner inequality and by using AI. 
We say a graph $H$ (as a subgraph of $G$) is a {\em cone} centered at $v$ if $v$ is connected to all other vertices of $H$. In such a case we write $H=\{v\}\cup H'$ and we call vertices of $H'$ \emph{leaves}. To extend Theorem \ref{cor:main-spectral} first we  analyzed cones which are ``stars'', i.e., the induced graph on $H'$ is an independent set. Observe that if $G$ is triangle free every cone is of this kind; this was also motivated by the fact that a large family of prior works studied triangle free graphs \cite{HV06}. The main advantage of studying a star is that in a typical proper coloring we expect many of the leaves to share colors and that leaves significantly more than $s$ colors to the center (so we should be able to implicitly improve upon the analysis of Lemma \ref{lem:mixed-form-lower-bound} above).

Then, we analyze the more general case of a cone. Roughly speaking, it seems we should ``lose" now at the center because, e.g., if $H'$ is a complete graph all leaves have distinct colors so the center may always only have $s$ colors left. But, in such a case we have many triangles involving the center and we can implicitly use a variant of Lemma~\ref{lem:mixed-form-lower-bound} to win.
Putting these together, for {\em any} cone and any function $f$ we obtain 
\begin{restatable}[Main Technical Theorem]{theorem}{localCorollary}\label{thm:local}
     For any $t>0$, and cone $H= \{ v \} \cup H'$ with center $v$ and leaves $H'$
    (edges among the leaves are allowed)
    with $|H'|=d$, suppose $s \geq 5$ and $|\cL_i| \geq s + \deg_H(i)$ for $i\in H$. Then we have
    \[
    2 \langle D_v f, D_{H'} f \rangle_{\mu}
    + t \|D_v f \|_{\mu}^2
    + \frac{\beta_{d}}{t(s-1)} \|D_{H'} f \|_{\mu}^2
    + O \Big( \frac{\Delta}{ts^3} \Big)
    \langle f, D_{H'} f \rangle_{\mu}
    \geq 0.
    \]
    where $D_{H'}=\sum_{i\in H'}D_i$ and $\beta_d=1-(\frac{s-1}{s})^d$.
\end{restatable}

Summing up the above inequality over a carefully chosen distribution of cones proves our main theorem. 
\begin{theorem}[Main]\label{thm:main}
    There exists an absolute constant $C>0$ such that, whenever $s\geq 5$ and $|\cL_i|\geq s+\deg(i)$ for all $i\in G$, we have 
    \[ D^2 \succeq \Bigg( 1 - \sqrt{\frac{\Delta \beta_{\Delta}}{s} + \frac{C\Delta^2}{s^3}} \Bigg) D =: (1-\gamma) D.\]
    Consequently, if $q>\alpha\Delta$ where $\alpha>\alpha^*$ for $\alpha^*
        =2-\exp\big(-1/(\alpha^*-1)\big)
        \approx 1.7408$ and sufficiently large $\Delta$, the spectral gap of the Glauber dynamics is at least $\tfrac{c_{\alpha}}{n}$ for some constant $c_{\alpha}>0$.
\end{theorem}

We emphasize that the above bound on $s$ is not optimal. We think of this work as a new method to approach the local-to-global method for graph coloring with the hope that our method leads to an optimal analysis of the Glauber dynamics. 

\subsection{Proof overview}
\label{sec:proofoverview}
In this section we explain the main ideas to prove Theorem~\ref{thm:local}. Fix an arbitrary cone with center $v$, and leaves $H'$. We assume the vertices outside $H$ are pinned to some partial proper coloring. So, in the worst case every leaf has at least $s+1$ colors available and the center $v$ has at least $s+|H'|$ many colors.

To prove Theorem \ref{thm:local} we mainly need to study the action of the operator $D_vD_{H'}$ on a generic function $f\in L^2(\mu)$. Our first observation is that it is enough to study the action on leaf functions; 
recall that $f$ is a leaf function if it is only a function of (the colors of) vertices in $H'$.
To do that, given $f$, we carefully choose a  leaf function $h^*$ (by solving an optimization problem) such that $P_v D_H (f-h^*)=0$. 
We further observe that $D_v h^*=0$ (this is because $D_v$ measures how $h^*$ changes when the color of the center changes, but $h^*$ is a leaf function so it is invariant under changes of the center). Having these, we get $\langle D_vf,D_{H'}f\rangle_\mu \geq \langle D_vf, D_vD_{H'}h^*\rangle_\mu.$
We then use a careful AM-GM to de-couple $D_vf$. Then, the main task reduces to upper-bounding $\norm{D_vD_{H'} h^*}_\mu.$

We are left with studying the action of the $D_vD_{H'}$ operator on a generic leaf function $h$. The first observation is that if $D_{H'} h$ is also a leaf function we are done. So, the whole argument is to say that $D_{H'} h$ is not ``far'' from the space of leaf functions. This is what we do in Lemma \ref{lem:approx-leaf}. In particular, we carefully choose a leaf function $h'$ and we show that $D_{H'}h$ is close to it.

First, we define a new operator $D_{H'}^{-v}$ which has no dependency on the center. In particular, $D_{H'}^{-v}h$ acts the same as $D_{H'}h$ except that we assume  the center is pinned to some ``dummy'' color that doesn't show up in any of the lists of the leaves. We now let $h'=D_{H'}^{-v}h$. Note that by definition $h'$ is a leaf function. 
 Now, the main observation is that for any color $c\in \cL_v$ we can write (see \eqref{eq:ell-DH})
$$(h'-D_{H'}h)(c,\bx)=(D_{H'}^{-v}h-D_Hh)(c,\bx)=\begin{cases}\sum_{i\in H'} \Pi_{i,c}^{-v}h(\bx)&\text{if }c\notin\bx\\0&\text{otherwise}\end{cases},$$ 
where $\Pi_{i,c}^{-v}$'s are rank-1 orthogonal projections onto the functions $e_{i,c}^{-v}$. The latter is a rank 1 centered function (i.e., expectation 0) that assigns color $c$ to vertex $i\in H'$ (assuming the center has a dummy color) and if $c\notin \cL_i$, $\Pi_{i,c}^{-v}=0.$
Having this, the problem essentially reduces to upper-bounding
\begin{equation}\label{eq:cxpiv}\sum_{c\in \cL_v} \E_{\bx}\bigg[ \mathbb{I}(c\notin \bx)\Big(\sum_{i\in H'} \Pi_{i,c}^{-v }h \Big)^2 \bigg] _\mu.\end{equation}
Observe that the inner-sum is just an expression on the subgraph $H'$; so we have in effect gotten rid of the center. 
The rank 1 operators, $\Pi_{i,c}$  give an ``orthonormal'' decomposition of the (normalized) Laplacians, $D_i$, namely $D_i\propto \sum_{c\in \cL_i} \Pi_{i,c}.$ (see Fact \ref{fact:sum-c}). To put things into perspective, we started with the goal of upper-bounding a term like $\E_{(c,\bx)}[(Dh)^2]$ over the graph $H$, and we have reduced to a similar quadratic quantity over $H'$, but only allowing for $\bx$'s that don't have $c$. In other words, instead of directly bounding Bochner type cross terms $\langle D_if,D_jf \rangle$, we need to bound $\langle \Pi_{i,c}f,\Pi_{j,c}f\rangle$ for vertices $i,j$ that share a common endpoint. We expect that this methodology will lead to a tight analysis of Glauber dynamics for graph coloring in future.  

We then upper-bound \eqref{eq:cxpiv} in Lemma \ref{lem:ABF}. The main observation is the following: Since every vertex $i\in H'$ has at least $s+\deg_{H'}(i)$ colors available, the probability that a uniformly random proper coloring of $H'$ doesn't assign the color $c$ to any vertex is at least $\Big(\frac{s-1}{s}\Big)^{|H'|}.$
The lemma proves the upper-bound 
$\beta_{|H'|}\bigg\| \sum_i \Pi^{-v}_{i,c}h \bigg\|^2$ on \eqref{eq:cxpiv} by writing the uniform distribution over proper colorings as a noise operator applied to the uniform distribution over proper colorings that don't use color $c$. In particular, to sample a proper coloring $\by$, we  first sample a proper coloring $\bx$ conditioned on $c\notin \bx$; then, we independently let $y_i=c$ with probability $\frac{1}{|\cL_i^{-v}|}$ and leave  $y_i=x_i$ otherwise. This allows us to use Cauchy-Schwarz to prove the lemma. One last technical property is that since $\Pi^{-v}_{i,c}$ operators are centered, we can essentially replace the event $\P_x[c\notin \bx]$ with $1-\P_x[c\notin \bx]$. 

Lastly, to finish the proof of Theorem \ref{thm:local} we need to prove an upper-bound on $\norm{\sum_{i\in H'} \Pi^{-v}_{i,c} h}^2.$ This is  done in Section \ref{sec:pairlocalineq}. In particular, the crux of the proof is Lemma \ref{lem:edge-frame} in which we upper-bound a cross term $\sum_c\langle \Pi_{i,c}f,\Pi_{j,c}f\rangle$ by $\langle D_1f,D_2f\rangle$ up to some lower order error terms. Although the proof is rather involved it builds on analysis of links of co-dimension 2 in the coloring complex as was done in the prior works \cite{ALO21}. The proof follows by decomposing $f$ as a sum of the degree 1 part and higher degree part. It turns out that the worst case is when $f$ is only a degree one function, i.e., $f(x_i,x_j)=\sum_c \alpha_{i,c}\mathbb{I}[x_i=c]+\sum_c \beta_c\mathbb{I}[x_j=c].$

\subsection{Concurrent work}
Independent of our work, Chen and Liu proved a near-optimal bound on the mixing time of the Glauber dynamics for graphs of girth $>4$ \cite{CL26a}. In another line, Haqi \cite{Haq26} independently proved that the Glauber dynamics mixes rapidly to generate a proper coloring of line graphs if $q\geq 2\Delta+4$.

\subsection{Acknowledgements}

The proof of Theorem~\ref{cor:main-spectral} is human made. The proof of Theorem \ref{thm:local} follows from a long interaction between human and AI. We started by first understanding an improved bound on sampling proper colorings of triangle free graphs and then we tried to extend the argument to all graphs.

The research is supported by a Simons Investigator Award 928589, an AFOSR grant FA9550261B166, and a Lazowska Endowed Professorship in Computer Science \& Engineering. AL is also supported by the Computer Science \& Engineering Educators' Endowed Fellowship and the Jonathan Hsieh and Kathleen Ting Endowed Fund. 

\section{Preliminaries}

We define $D_i = I - P_i$ where $P_if = \E[f\, | \, \bX_{-i}]$. Here $P_i,D_i \colon L^2(\mu) \to L^2(\mu)$. 

For a matrix $M\in \R^{m\times n}$ we write $\norm{M}:=\max_{\bx\in \R^n \smallsetminus \{0 \}}\frac{\norm{M\bx}}{\norm{\bx}}.$ We also recall $\norm{M}_1$ is the maximum sum of absolute values of entries of a row of $M$ and $\norm{M}_\infty$ is the maximum sum of absolute values of entries of a column. For such a $M$, one can write $\norm{M}\leq\sqrt{\norm{M}_1\norm{M}_\infty}$.

For a vertex $i$ and a subgraph $H \subset G$ we write $N_H(i)$ for the neighbors of $i$ in $H$.
 
Given a graph $G$ where every vertex $i$ has a list of colors $\cL_i$, for a set of colors $C$ we write $\cL_i^{-C}$ to denote $\cL_i\smallsetminus C$.  Similarly, for a vector $\bx$ we write $\cL_i^{-\bx_{N(i)}}$ to denote $\cL_i\smallsetminus \bx_{N(i)}$. 
Analogously, we write $\ell_i$ or $\ell_i^{-C}$ to denote the sizes of $\cL_i$, $\cL_i^{-C}$ respectively.

For a set $S$, we write $X \sim S$ to denote that  $X$ is chosen to be a uniformly random element of $S$.

\begin{lemma}[Matrix Hoeffding Inequality]\label{lemma:Hoeffding}
Let $G$ be a graph and $\cH$ be a real Hilbert space, and suppose
\[
z=(z_v)_{v\in V(G)}\in\cH^{V(G)}.
\]
Then we have
\[
\sum_{i,j\in G}A_G(i,j)\langle z_i,z_j\rangle_{\cH}
\ge
\lambda_{\min}(A_G)\sum_{i\in G}\Norm{z_i}_{\cH}^2.
\]
\end{lemma}

\subsection{Local Inequalities over Subgraphs}
Throughout the paper we study local inequalities for a subgraph  $H \subset G$; in doing so, we always assume that we have fixed a partial proper coloring for vertices in $G\smallsetminus H$. Thus, each vertex $i\in H$ is left with a list $\cL_i$ of colors (after removing the conflicts with outside vertices). We define 
\[ \Omega_H = \bigg\{\bx \in \prod_i \cL_i: x_i \neq x_j \; \forall ij \in E(H)  \bigg\} \quad\text{and}\quad \mu_H = \Unif(\Omega_H).\]

We note that operators $D_i,P_i$ naturally localize to this subgraph:
 \[(P_i^Hf)(\bx) = \E_{\mu_H}[f \, | \, \bX_{-i} = \bx_{-i}] = \frac{1}{\ell_i^{-\bx_{N_H(i)}}} \sum_{c \in \cL_i^{-\bx_{N_H(i)}}}f(\bx_{i \to c}) \quad\text{and} \quad (D_i^Hf)(\bx) = f(\bx) - (P_i^Hf)(\bx), \]
 so we can just remove the subscript $H$. 
 We next discuss a decomposition statement of $D_i$ into
rank-one projection matrices.

Conditioned on $\bX_{-i}=\bx_{-i}$,  observe that the  distribution of $X_i$
is uniform on $\cL_i^{-\bx_{N_{H}(i)}}$.
For $c \in \cL_i^{-\bx_{N_{H}(i)}}$, we define $p_{i,c}(\bx)=\P_\mu[X_i=c\mid\bX_{-i}=\bx_{-i}]$. 
Let
\[
e_{i,c}(\bx)
=\frac{\mathbb{I}\{x_i=c\}-p_{i,c}}
       {\sqrt{p_{i,c}(1-p_{i,c})}}
=\begin{cases}
 \sqrt{\dfrac{1-p_{i,c}}{p_{i,c}}},&x_i=c,\\[1ex]
 -\sqrt{\dfrac{p_{i,c}}{1-p_{i,c}}},&x_i\ne c.
\end{cases}
\]
Otherwise, if $c  \not \in \cL_i^{-\bx_{N_{H}(i)}}$, we define $e_{i,c}=0$. See section \ref{sec:proofoverview} for motivation. Also, notice that the vector $e_{i,c}$ naturally localizes even to subgraphs of $H'\subset H$ (once we fix a partial proper coloring of $H\smallsetminus H'$).

For $c \in \cL_i^{-\bx_{N_H(i)}}$, we note that
\[
\E_\mu[e_{i,c}(\bX)\mid\bX_{-i}=\bx_{-i}]=0 \quad \text{and} \quad
\E_\mu[e_{i,c}(\bX)^2\mid\bX_{-i}=\bx_{-i}]=1.
\]
 
We naturally assume all these functions are defined over $H$ and throughout the paper we may localize them further to subgraphs of $H$.

For a vertex $i$, and $H=\{i\}$ note that the functions $\mathbb{I}\{x_i = c\}$'s for $c\in \cL^H_i$ are orthogonal under $\mu$. But because $D^H_i$ zeroes out the all-ones function we also normalize and project out the all-ones vector from each indicator function.
The resulting vectors $e_{i,c}$'s don't form an orthonormal basis (because they form $\ell_i$ vectors in an $\ell_i-1$ dimensional space) but they span the space of functions orthogonal to all ones. 

More generally, for $H$ arbitrary and $i\in H$, the following properties are immediate: 
\begin{itemize}
\item For distinct colors $c,c'$ that are both available,
we have
\begin{equation}\label{eq:eic-inner}
\E_\mu[e_{i,c}(\bX)e_{i,c'}(\bX)\mid\bX_{-i}=\bx_{-i}]
=-\sqrt{\frac{p_{i,c}p_{i,c'}}
              {(1-p_{i,c})(1-p_{i,c'})}}
=-\frac1{\ell_i^{-\bx_{N_H(i)}}-1}.
\end{equation}
\item For any $c\in\cL_i^{-\bx_{N_H(i)}}$,
\[
D_i e_{i,c}=e_{i,c}.
\]
\end{itemize}

Now, we let $\Pi_{i,c}$ be the orthogonal projection onto $e_{i,c}$, when $c \in \cL_i^{-\bx_{N_H(i)}}$ and $0$ otherwise. 
In particular, for any $h\in L^2(\mu)$
\begin{equation}\label{eq:prodi}
(\Pi_{i,c}h)(\bx)
:= \begin{cases} e_{i,c}(\bx)\,
  \E_\mu [h(\bX)e_{i,c}(\bX)\mid\bX_{-i}=\bx_{-i}] &\text{if } c \in \cL_i^{-\bx_{N_{H}(i)}}, \\ 0 &\text{otherwise.} \end{cases}
\end{equation}

We note that we might sometimes drop the superscript of $H$ when the corresponding subgraph is clear from context. 

For a graph consisting of one edge on vertices $1$ and $2$, we fix $c \in \cL_1 \cap \cL_2$. For every proper
coloring $(x,y) \in \Omega$ with $x,y \neq c$, we have (recall that $e_{1,c}(x,y) = \tfrac{\ell_1^{-y}\mathbb{I} \{x = c \} - 1}{\sqrt{\ell_1^{-y}-1}}$)
        \begin{align}\label{eq:explicit-Pi-i}
            (\Pi_{1,c} f)(x,y) &= \E_{\mu}[f e_{1,c} \, | \, X_2 = y] e_{1,c}(x,y) \notag \\
            &=\frac{f(c,y) - (P_1f)(y)}{\sqrt{\ell_1^{-y} - 1}} \cdot \frac{-1}{\sqrt{\ell_1^{-y}-1}} =- \frac{D_1f(c,y)}{\ell_1^{-y} - 1}.
        \end{align}
Furthermore, by definition if $c= y$ (resp. $c=x$) then $\Pi_{1,c}f(x,c) = 0$ (resp. $\Pi_{2,c} f(c,y) = 0$). 

\begin{fact}\label{fact:sum-c}
For any $i$, any $h \in L^2(\mu)$, and any $\bx$, we have
\[
\sum_{c\in\cL_i}(\Pi_{i,c}h)(\bx)
=\left(1+
 \frac1{\ell_i^{-\bx_{N_{H}(i)}}-1}\right)
 (D_i h)(\bx).
\]
If $H$ is an independent set, $\ell^{-\bx_{N_H(i)}}=\ell^H_i$ and we get
$\sum_{c\in\cL_i}\Pi_{i,c}
=\frac{\ell_i}{\ell_i-1}D_i.$
\end{fact}

\begin{proof}
We fix $\bx_{-i}$ and let $c' \in \cL_i^{-\bx_{N_H(i)}}$. We may write
\begin{align*}
\left(\sum_{c\in\cL_i}\Pi_{i,c}\right)e_{i,c'}
&=\sum_{c\in\cL_i}
 \E_\mu[e_{i,c}(\bX)e_{i,c'}(\bX)
        \mid\bX_{-i}=\bx_{-i}]\,e_{i,c}\\
&\underset{\eqref{eq:eic-inner}}{=}
 e_{i,c'}
 -\frac{
        \sum_{c\in\cL_i, c\ne c'}e_{i,c}}
       {\ell_i^{-\bx_{N_H(i)}}-1}\\
&=\left(1+
 \frac1{\ell_i^{-\bx_{N_H(i)}}-1}\right)e_{i,c'}.
\end{align*}
Since the vectors $e_{i,c'}$ for the available colors span the space
of functions orthogonal to $\mathbf{1}$, and $\Pi_{i,c}$ and $D_i$ annihilate $\mathbf{1}$, the desired conclusion
follows.
\end{proof}

\section{Proof of Theorem~\ref{thm:main}}
In this section we prove Theorem \ref{thm:main}. As we said in the introduction, the proof is based on a careful averaging argument over cones. 

\begin{proof}[Proof of Theorem~\ref{thm:main} assuming Theorem~\ref{thm:local}]
    We assume that $\gamma <1 $ as otherwise Theorem~\ref{thm:main} is vacuously true. The plan is to average Theorem~\ref{thm:local} over random cones and sum them up. First, we choose ${\bf p} \colon E(G) \to  [1- \kappa/2,1]$  such that for every edge $e=uv$, 
    \begin{equation}\label{eq:puv}
        p_{uv}+ \frac{\kappa}{2\Delta}\sum_{w \in N_G(u) \cap N_G(v)} p_{wu}p_{wv} = 1,
    \end{equation}
    where $\kappa$ will be defined later such that $\kappa/2<1$.
    Note that by definition of ${\bf p}$, $p_{uv}=p_{vu}$. We defer the justification for the existence of ${\bf p}$ to the end of the proof. 
    For every vertex $v$, we choose the corresponding cone $H_v := G[\{ v\} \cup H_v' ]$ to apply Theorem~\ref{thm:local} as follows: each $u \in N(v)$ is included in $H_v'$ {\em independently} with probability $p_{uv}$. 
    
    By  Theorem~\ref{thm:local} with $t = \gamma$, there exists some $C'>0$ such that for every $f \colon \Omega \to \R$, 
    \begin{equation}\label{eq:sum-local}
        \frac{1}{2} \sum_v \E\bigg[2 \langle D_vf, D_{H_v'}f \rangle_{\mu} + \gamma \| D_v f\|_{\mu}^2 + \frac{\kappa}{ \Delta} \|D_{H_v'} f \|_{\mu}^2 + \frac{C'\Delta}{\gamma s^3} \langle f, D_{H_v'} f \rangle\bigg] \geq 0,
    \end{equation}
    where we write 
    $\kappa = \tfrac{\Delta \beta_{\Delta}}{\gamma(s-1)}$ as a shorthand. Note that $0 \leq  \kappa \leq \gamma < 1$.

    Since $D_v^2 = D_v$, the coefficient of $D_v$ in the LHS of \eqref{eq:sum-local} is 
    \begin{equation}\label{eq:coeff-Dv}
        \frac{\gamma}{2} + \frac{1}{2}\bigg(\frac{\kappa}{\Delta} + \frac{C' \Delta}{\gamma s^3} \bigg) \sum_{u \in N_G(v)} p_{uv} \underset{\sum_{u \in N_G(v)} p_{uv}  \leq \Delta}{\leq} \frac{\gamma+\kappa + \frac{C' \Delta^2}{\gamma s^3}}{2} \leq \gamma,
    \end{equation}
    where the final inequality follows by choosing for instance $C \geq C'+2$. 
    
    For an edge $uv$, the coefficient of $\langle D_uf, D_vf \rangle$ in the LHS of \eqref{eq:sum-local} is 
    \begin{equation}\label{eq:coeff-DuDv}
     (p_{uv}+p_{vu})+
\frac {\kappa}{\Delta}\sum_{w\in N_G(u)\cap N_G(v)}p_{wu}p_{wv} \underset{\eqref{eq:puv}, \, p_{uv} = p_{vu}}{=}2.
    \end{equation}
    For non-edges $uv$,  $D_u$ and $D_v$ are commuting orthogonal
projections so that for any $f \colon \Omega \to \R$ we have $\langle D_uf, D_v f \rangle \geq 0$ and furthermore the term $\langle D_uf, D_v f \rangle$ appears in the LHS of \eqref{eq:sum-local} with coefficient 
\begin{equation}\label{eq:coeff-nonedge}
    \frac{\kappa}{\Delta}
\sum_{w\in N_G(u)\cap N_G(v)}p_{wu}p_{wv} \in [0, \kappa] \subset [0,2],
\end{equation}
where we use the fact that $p_{uw} \in [0,1]$ for $w \in N_G(u)$ and $|N_G(u) \cap N_G(v)| \leq \Delta$. 

Consequently, \eqref{eq:sum-local} in turn implies that for any $f \colon \Omega \to \R$, we have 
\[\langle f, (D^2 - (1- \gamma)D) f \rangle_{\mu}=\gamma \sum_v \| D_v f\|_{\mu}^2 + 2\sum_{ u \neq v \text{ unordered}} \langle D_uf, D_v f\rangle_{\mu}   \underset{\eqref{eq:coeff-Dv}, \, \eqref{eq:coeff-DuDv}, \,\eqref{eq:coeff-nonedge}}{\geq} 0.\]
Finally, we show that there exists ${\bf p}$ satisfying the property in \eqref{eq:puv}. To that end, we define a  continuous map
$T:[1-\kappa/2,1]^{E(G)}\to [1-\kappa/2,1]^{E(G)}$. For a vector ${\bf p}$ we write,
\[
(T{\bf p})_{uv}:=1-\frac{\kappa}{2 \Delta}
\sum_{w\in N_G(u)\cap N_G(v)}p_{wu}p_{wv} \in 
[1 - \kappa/2, 1],\quad\forall uv\in E(G).
\]
where we used that $u,v$ have at most $\Delta$ common neighbors. Since $[1-\kappa/2,1]^{E(G)}$ is a nonempty compact convex set, Brouwer's fixed point theorem gives a
fixed point ${\mathbf p} \in [1-\kappa/2,1]^{E(G)}$ satisfying the desired conditions.

Finally, to deduce the spectral gap of Glauber dynamics when $q > \alpha \Delta$ for $\alpha > \alpha^*$, we apply the bound above with  $s=(\alpha-1)\Delta$. We note that $\gamma^2
=\frac{1-\left(1-\frac{1}{(\alpha-1)\Delta}\right)^\Delta}
        {\alpha-1}
  +\frac{C}{(\alpha-1)^3\Delta} \xrightarrow[\Delta\to\infty]{}
  \frac{1-e^{-1/(\alpha-1)}}{\alpha-1}<1 $. The final inequality follows because $1-e^{-1/x}$ is decreasing for $x>0$ and equals
$\alpha^*-1$ at $x=\alpha^*-1$. So for sufficiently large $\Delta$,
we have $1-\gamma\ge c_\alpha>0$. By Fact~\ref{fact:LLvsL}, it follows in this case that the spectral gap of Glauber dynamics $\ge \tfrac{c_{\alpha}}{n}$, as desired.
\end{proof}

\section{Pairwise Local Inequalities}\label{sec:pairlocalineq}
The following is the main statement of this section. We prove it using Lemmas \ref{lem:non-edge-frame} and \ref{lem:edge-frame} proven below.

\begin{lemma}\label{lem:projection-sum}
If $\ell_i \geq s + \deg_H(i)$ for each $i \in H$ and $s \geq 5$, then for every $h\in L^2(\mu_H)$, we have
\[
\sum_c\left\|\sum_{i\in H}\Pi_{i,c}h\right\|_{\mu_H}^2
\le\left(1+\frac1{s-1}\right)\|Dh\|_{\mu_H}^2
+\frac{O(1)}{s^2}
 \sum_{i\in H}\deg_H(i)\|D_ih\|_{\mu_H}^2.
\]
\end{lemma}

\begin{proof}
Since each $\Pi_{i,c}^H$ is an orthogonal projection, Fact~\ref{fact:sum-c} gives
\[
\sum_c\|\Pi_{i,c}h\|_{\mu_H}^2
=\bigg\langle h,\sum_c\Pi_{i,c}h\bigg\rangle_{\mu_H}
\underset{\ell_i^{-x_{N_H(i)}}\ge s}{\leq}\left(1+\frac1{s-1}\right)\|D_ih\|_{\mu_H}^2.
\]

Fix distinct nonadjacent vertices $i,j\in H$. After pinning the other
vertices, both remaining lists have size at least $s$. Consequently, by averaging over the pinning,
we obtain
\[
\sum_c\langle\Pi_{i,c}h,\Pi_{j,c}h\rangle_{\mu_H}
\underset{\text{Lemma~\ref{lem:non-edge-frame}}}{\le}\left(1+\frac1{s-1}\right)
      \langle D_ih,D_jh\rangle_{\mu_H}.
\]

Similarly for adjacent vertices $i,j\in H$, after pinning the other vertices, both remaining lists have size at least $s+1$, and so by averaging over the pinning, we obtain
\begin{align*}
\sum_c\langle\Pi_{i,c}h,\Pi_{j,c}h\rangle_{\mu_H}
&\underset{\text{Lemma~\ref{lem:edge-frame}}}{\le}\langle D_ih,D_jh\rangle_{\mu_H}
 +\frac{O(1)}{s^2}
  \bigl(\|D_ih\|_{\mu_H}^2+\|D_jh\|_{\mu_H}^2\bigr)\\
&\underset{\text{Lemma~\ref{lem:mixed-form-lower-bound}}}{\le}\left(1+\frac1{s-1}\right)
      \langle D_ih,D_jh\rangle_{\mu_H}
 +\frac{O(1)}{s^2}
  \bigl(\|D_ih\|_{\mu_H}^2+\|D_jh\|_{\mu_H}^2\bigr).
\end{align*}
Summing all these proves the lemma.
\end{proof}

\begin{lemma}\label{lem:non-edge-frame}
Given a graph of two isolated vertices labeled 1 and 2, suppose $s\geq 2$, and let $\mathcal L_1,\mathcal L_2$
be lists of sizes $\ell_1,\ell_2\geq s$. Then for every $f \in L^2(\mu)$, we have
\begin{equation}\label{eq:non-edge-frame}
\sum_c\langle\Pi_{1,c}f,\Pi_{2,c}f\rangle_\mu
\leq
\left(1+\frac{1}{s-1}\right)\langle D_1f,D_2f\rangle_\mu.
\end{equation}
\end{lemma}

\begin{proof}
The difference $D_2-\Pi_{2,c}$ is
an orthogonal projection since $\Pi_{2,c}$ is the projection onto $e_{2,c}$ and $P_2e_{2,c} = 0$. Similarly, $D_1 - \Pi_{1,c}$ is an orthogonal projection.

    Furthermore, $D_1D_2 \succeq 0$, $\Pi_{1,c}(D_2-\Pi_{2,c}) \succeq 0$  and $\Pi_{2,c}(D_1-\Pi_{1,c}) \succeq 0$  because operators acting
on different coordinates commute under the product measure. This allows us to write
\begin{align*}
    &\left(1+\frac1{s-1}\right)D_1D_2
      -\sum_c\Pi_{1,c}\Pi_{2,c}\\
    &\underset{\text{Fact~\ref{fact:sum-c}}}{=}
      \left(1+\frac1{s-1}-\frac{\ell_1}{2(\ell_1-1)} - \frac{\ell_2}{2(\ell_2-1)}\right)
      D_1D_2
      +\frac{1}{2} \sum_c \bigg( \Pi_{1,c}(D_2-\Pi_{2,c}) + \Pi_{2,c}(D_1 - \Pi_{1,c})\bigg)\\
    &\quad=
      \left(\frac1{s-1}-\frac1{2(\ell_1-1)} - \frac{1}{2(\ell_2-1)}\right)
      D_1D_2
      +\frac{1}{2} \sum_c \bigg( \Pi_{1,c}(D_2-\Pi_{2,c}) + \Pi_{2,c}(D_1 - \Pi_{1,c})\bigg)
      \succeq_\mu0.
\end{align*}
In the first equality we are using the fact that, in this subgraph with two disconnected vertices, Fact~\ref{fact:sum-c} gives $\sum_c \Pi_{1,c}=\frac{\ell_1}{\ell_1-1}D_1$.
\end{proof}

\begin{restatable}{lemma}{edgeFrameLemma}
\label{lem:edge-frame}
Given a graph consisting of an \emph{edge} between vertices 1 and 2. Suppose $s\geq 5$, and let $\mathcal L_1,\mathcal L_2$
be lists of sizes $\ell_1,\ell_2\geq s + 1$. Then for every $f \in L^2(\mu)$, we have

\begin{equation}\label{eq:edge-frame}
\sum_c\langle\Pi_{1,c}f,\Pi_{2,c}f\rangle_\mu
\leq\langle D_1f,D_2f\rangle_\mu
+\frac{O(1)}{s^2}
 \big(\norm{D_1f}_\mu^2+\norm{D_2f}_\mu^2\big).
\end{equation}
\end{restatable}

\begin{proof}[Proof of Lemma~\ref{lem:edge-frame}]
    We define $\cS$ to be the restriction to $\Omega$ of functions $f_0 + f_1(X_1) + f_2(X_2)$ where $f_0 \in \R$ and $f_i \colon \cL_i \to \R$. We decompose
    \[ f = f_S + f_S^{\perp}, \quad \text{where } f_S \in \cS \text{ and } f_S^{\perp} \perp_{L^2(\mu)} \cS.\]
    Here, for any $(x_1,x_2)\in\Omega$ we may write, 
    \begin{equation}\label{eq:f1c0}
        f_S(x_1,x_2) = f_0 + f_1(x_1) + f_2(x_2), \quad\text{where } \sum_{c \in \cL_1} f_1(c) = \sum_{c \in \cL_2}f_2(c) = 0.
    \end{equation}

    We define an operator $\sol{P}_1$ that maps $L^2(\mu)\to\R^{\cL_2}$ where $\sol{P}_1g(c)=\E_{\mu}[g(X_1,X_2)|X_2=c]$. Similarly, we think about $\sol{P}_2$.
    
    We claim that 
    \begin{equation}\label{eq:P2D1-P1D2}
        \bigg| \sum_c \langle \Pi_{1,c} f, \Pi_{2,c} f \rangle_{\mu} - \frac{1}{|\Omega|} \sum_{ c \in \cL_1 \cap \cL_2} (\sol{P}_2D_1f)(c) \cdot (\sol{P}_1D_2f)(c) \bigg| \leq \frac{O(1)}{s^2} \|D_1 f \|_{\mu} \|D_2 f \|_{\mu}
    \end{equation}
    and 
    \begin{equation}\label{eq:D1D2-P2D1-P1D2}
        \bigg| \langle D_1 f_S, D_2f_S \rangle_{\mu} - \frac{1}{|\Omega|} \sum_{c \in \cL_1 \cap \cL_2}(\sol{P}_2D_1f)(c)\cdot (\sol{P}_1D_2f)(c) \bigg| \leq \frac{O(1)}{s^2} \|D_1 f_S \|_{\mu} \|D_2 f_S \|_{\mu}. 
    \end{equation}

    Assuming \eqref{eq:P2D1-P1D2} and \eqref{eq:D1D2-P2D1-P1D2}, we now complete the proof. Since $\cS$ contains every function of $X_1$, for every
$x\in\cL_1$ we have
\begin{equation}\label{eq:perp=0}
    (P_2f_S^{\perp})(\bx) = \frac{1}{\ell_2^{-x}} \cdot \sum_{y \in \cL_2 \setminus \{x \}} f_S^{\perp}(x,y) = \frac{|\Omega|}{\ell_2^{-x}} \cdot  \langle f_S^{\perp}, \mathbb{I}\{X_1 = x \} \rangle_{\mu} \underset{f_S^{\perp} \perp_{L^2(\mu)}\cS}{=} 0.
\end{equation}
Similarly, we have $P_1f_{S}^{\perp} = 0$. This implies that $D_1 f_S^{\perp} = D_2f_S^{\perp} = f_S^{\perp}$. Consequently, we have 
\begin{equation}\label{eq:quad-decompose}
    \langle D_1f,D_2f\rangle_\mu
    =\langle D_1f_S,D_2f_S\rangle_\mu
      +\norm{f_S^\perp}_\mu^2. 
\end{equation}

With this decomposition, it follows that 

\begin{align*}
&\sum_c\langle\Pi_{1,c}f,\Pi_{2,c}f\rangle_\mu - \langle D_1f,D_2f\rangle_\mu\\
    &\underset{\eqref{eq:quad-decompose}}{=}
      \bigg(
        \sum_c\langle\Pi_{1,c}f,\Pi_{2,c}f\rangle_\mu
        -\frac1{|\Omega|}
         \sum_{c\in\cL_1\cap\cL_2}(\sol{P}_2D_1f)(c)\cdot(\sol{P}_1D_2f)(c)
      \bigg) - \norm{f_S^\perp}_\mu^2\\
    &\qquad-
     \bigg(
        \langle D_1f_S,D_2f_S\rangle_\mu - \frac{1}{|\Omega|}\sum_{c\in\cL_1\cap\cL_2}(\sol{P}_2D_1f)(c)\cdot(\sol{P}_1D_2f)(c)\bigg) \\
    &\underset{\eqref{eq:P2D1-P1D2}, \, \eqref{eq:D1D2-P2D1-P1D2}}{\leq} \frac{O(1)}{s^2} \|D_1f\|_{\mu} \|D_2 f \|_{\mu} + \frac{O(1)}{s^2} \|D_1f_S\|_{\mu} \|D_2 f_S \|_{\mu} - \|f_S^{\perp} \|_{\mu}^2 \underset{\|D_if_S \|_{\mu} \leq \|D_if \|_{\mu}}{\leq} \frac{O(1)}{s^2} \|D_1f \|_{\mu}\|D_2 f \|_{\mu}. 
\end{align*}

\paragraph{Proof of \eqref{eq:P2D1-P1D2}} 
First, we note that 
\begin{equation}\label{eq:LHS-Pi1Pi2}\sum_c\langle\Pi_{1,c}f,\Pi_{2,c}f\rangle_\mu = \frac{1}{|\Omega|} \sum_{c \in \cL_1 \cap \cL_2} \sum_{(x,y) \in \Omega}  (\Pi_{1,c} f)(x,y) \cdot (\Pi_{2,c} f)(x,y) \underset{\eqref{eq:explicit-Pi-i}}{=} \frac{1}{|\Omega|} \sum_{c \in \cL_1 \cap \cL_2} \sum_{\substack{x \in \cL_1, \, y \in \cL_2 \\x,y,c \text{ distinct}}} \frac{D_1f(c,y) D_2f(x,c)}{( \ell_1^{-y} - 1 ) (\ell_2^{-x} - 1)}.
\end{equation}
where the second identity uses that $\Pi_{1,c}f(x,c)=\Pi_{2,c}f(c,y)=0$.
        
Having this we prove \eqref{eq:P2D1-P1D2}. For $c \in \cL_1 \cap \cL_2$, we note that
\[ (\sol{P}_2D_1f)(c) = \frac{1}{\ell_2 - 1 } \sum_{y \in \cL_2 \setminus \{c \}} D_1f(c,y) \quad\text{and}\quad (\sol{P}_1D_2f)(c) = \frac{1}{\ell_1 - 1}\sum_{x \in \cL_1 \setminus \{c \}}D_2f(x,c).\]
Let $\Omega_1,\Omega_2\subset \Omega$ be the sets of all proper colorings $(x,y)$ such that $x\in \cL_1\cap \cL_2$ (resp. $y\in\cL_1\cap\cL_2$). Let $g_1\in\R^{\Omega_1}$ be equal to $D_1f$ and $g_2\in \R^{\Omega_2}$ be $D_2f$. Then, 
\begin{align}\label{eq:M-P2D1}
    &\sum_{c}\langle \Pi_{1,c} f, \Pi_{2,c} f \rangle_{\mu} - \frac{1}{|\Omega|} \sum_{c \in \cL_1 \cap \cL_2} (\sol{P}_2D_1f)(c) \cdot (\sol{P}_1D_2f)(c) \notag \\
    &=
    \frac{1}{|\Omega|}\sum_{c\in \cL_1\cap \cL_2}\sum_{x\in \cL_1^{-c}, \, y\in \cL_2^{-c}}g_1(c,y)g_2(x,c)\left(\frac{\mathbb{I}\{x \neq y\}}{(\ell_1^{-y}-1)(\ell_2^{-x}-1)}-\frac{1}{(\ell_1-1)(\ell_2-1)}\right)\\
    &\leq\sum_{c\in \cL_1\cap \cL_2}\frac{\norm{M^c}}{|\Omega|} \cdot\norm{g^c_1}\norm{g^c_2},
\end{align}
where the norms are the Euclidean $\ell_2$-norm, $g_1^c, g_2^c$ is the restriction of $g_1,g_2$ to proper colorings where $1$ (resp $2$) is colored $c$. Furthermore, for each $c\in \cL_1\cap \cL_2$, define $M^c \in \R^{\cL_2^{-c} \times \cL_1^{-c}}$ by 
\[ (M^c)_{y,x} = \frac{\mathbb{I} \{x \neq y \}}{(\ell_1^{-y} -1)(\ell_2^{-x} - 1)} - \frac{1}{(\ell_1 - 1)(\ell_2-1)}.\]

Fix $c\in \cL_1\cap \cL_2$. Then, $\norm{M^c}\leq \sqrt{\norm{M^c}_\infty\norm{M^c}_1}$. By definition, $\|M^c\|_1 =\max_{y\in \cL_2^{-c}}\sum_{x \in \cL_1^{-c}}|M^c_{y,x} | $. For $y\in \cL_2^{-c}$ arbitrary we have
\begin{align*}
     \sum_{x \in \cL_1^{-c}}|M^c_{y,x} | &= \sum_{x \in \cL_1^{-c-y}} \bigg( \frac{1}{(\ell_1^{-y} - 1)(\ell_2^{-x} - 1)} - \frac{1}{(\ell_1 - 1)(\ell_2-1)}\bigg) + \frac{\mathbb{I}\{y \in \cL_1 \}}{(\ell_1 - 1)(\ell_2-1)} \leq \frac{O(1)}{s^2}
\end{align*}
where the inequality uses for $x \neq y$ that $0\le M^c_{y,x}
\le
\frac{s}{\ell_1-1}
\left(\frac{1}{(s-1)^2}-\frac{1}{s^2}\right)
\le
\frac{O(1)}{(\ell_1-1)s^2}$.
 Analogously,  we have $\|M^c \|_{\infty} \leq \frac{O(1)}{s^2}$. Therefore, $\| M^c \|\leq \frac{O(1)}{s^2}$.

Consequently, we have 
\begin{align*}
    \bigg|  \sum_{c}\langle \Pi_{1,c} f, \Pi_{2,c} f \rangle_{\mu} - \frac{1}{|\Omega|}& \sum_{c \in \cL_1 \cap \cL_2} (\sol{P}_2D_1f)(c)(\sol{P}_1D_2f)(c)  \bigg| \underset{\eqref{eq:M-P2D1}}{\leq} \sum_{c\in\cL_1\cap \cL_2}\frac{\| M^c \|}{|\Omega|} \norm{g_1^c}\norm{g_2^c} \\
    &\underset{\text{Cauchy-Schwarz}}{\leq} \frac{O(1)\norm{g_1}\norm{g_2} }{|\Omega|s^2}\leq \frac{O(1)}{s^2}\norm{D_1f}_\mu\norm{D_2f}_{\mu}. 
\end{align*}
This proves \eqref{eq:P2D1-P1D2}.

\paragraph{Proof of \eqref{eq:D1D2-P2D1-P1D2}}
We first claim that 
\begin{equation}\label{eq:P2D1-rearrange}
      \frac{1}{|\Omega|} \sum_{c \in \cL_1 \cap \cL_2} (\sol{P}_2D_1f)(c) \cdot (\sol{P}_1D_2f)(c) - \langle D_1 f_S, D_2f_S\rangle_{\mu} = \langle P_1f_1, P_2f_2 \rangle_{\mu} + \frac{1}{|\Omega|}\sum_{c \in \cL_1 \cap \cL_2} (\sol{P}_2P_1f_1)(c) \cdot (\sol{P}_1 P_2f_2)(c).
\end{equation}

Having this, and using $\PP_{\mu}[X_1 = c] = \frac{\ell_2 - 1}{|\Omega|}$ (and similarly swapping the roles of 1 and 2), we may write 
\begin{align}\label{eq:intermediate-P1P2f}
    &\bigg| \langle D_1f_S, D_2 f_S \rangle_{\mu} - \frac{1}{|\Omega|} \sum_{ c \in \cL_1 \cap \cL_2} (\sol{P}_2D_1f)(c) \cdot (\sol{P}_1D_2f)(c) \bigg| \notag \\&\hspace{-4em}\underset{\text{Triangle inequality}}{\leq} |\langle P_1f_1, P_2f_2 \rangle_{\mu}|+\frac{1}{\sqrt{(\ell_1 -1) (\ell_2-1)}} \bigg| \sum_{c \in \cL_1 \cap \cL_2} \sqrt{\PP_{\mu}[X_1 = c]}(\sol{P}_2P_1f_1)(c) \cdot \sqrt{\PP_{\mu}[X_2 = c]} (\sol{P}_1P_2f_2)(c) \bigg| \notag \\
    &\hspace{-2em}\underset{\text{Cauchy-Schwarz}}{\leq} \|P_1f_1\|_{\mu} \cdot \|P_2f_2\|_{\mu}+ \frac{1}{\sqrt{(\ell_1 - 1)(\ell_2-1)}}\bigg(\sum_{c \in \cL_1 \cap \cL_2} \P_{\mu}[X_1 = c] (\sol{P}_2 P_1f_1)(c)^2\bigg)^{1/2} \\ &\hspace{24em} \cdot\bigg( \sum_{c \in \cL_1 \cap \cL_2} \P_{\mu}[X_2 = c] (\sol{P}_1 P_2f_2)(c)^2 \bigg)^{1/2}\notag\\
    &\leq\|P_1f_1\|_{\mu} \cdot \|P_2f_2\|_{\mu}+\|P_2P_1f_1\|_{\mu} \cdot \|P_1P_2f_2\|_{\mu}\leq O(1)\cdot \|P_1f_1\|_{\mu}\cdot\|P_2f_2\|_{\mu}\\
    &\underset{\text{Claim~\ref{claim:P2f2-D2f2}}, \, \|D_2f_S \|_{\mu} = \|D_2 f_2 \|_{\mu}}{\leq} \frac{O(1) \cdot \|D_1 f_S \|_{\mu} \cdot \|D_2 f_S \|_{\mu}}{(\ell_1 - 1)(\ell_2-1)-1} \leq \frac{O(1)}{s^2} \|D_1 f_S \|_{\mu} \|D_2 f_S \|_{\mu}.
\end{align}
where we used that $P_2P_1f$ naturally extends $\sol{P}_2P_1f$ to $L^2(\mu)$.

It remains to prove \eqref{eq:P2D1-rearrange}. First, note that 
\begin{align}\label{eq:LHS-P2D1-identity}
    \langle D_1f_S, D_2f_S\rangle_{\mu} &= \langle f_1,f_2 \rangle_{\mu} - \langle f_1, P_1f_2 \rangle_{\mu} - \langle P_2f_1, f_2 \rangle_{\mu} + \langle P_1f_1, P_2f_2\rangle_{\mu} \notag \\
    &\underset{P_1f_2 = f_2, \, P_2f_1 = f_1}{=} - \langle f_1, f_2 \rangle_{\mu} + \langle P_1f_1, P_2 f_2 \rangle_{\mu} \notag \\
    &\underset{\text{Claim~\ref{claim:ip-sum}}}{=} \frac{1}{|\Omega|} \sum_{c \in \cL_1 \cap \cL_2} f_1(c)f_2(c) + \langle P_1f_1, P_2f_2 \rangle_{\mu}.
\end{align}

Next, we note that 
\[ D_1f = D_1f_S + D_1f_{S}^{\perp} \underset{\eqref{eq:perp=0}}{=} f_1 - P_1f_1 + f_S^{\perp}\]
which in turn shows that $\sol{P}_2D_1f = \sol{P}_2f_1 - \sol{P}_2P_1f_1 + \sol{P}_2f_S^{\perp} = f_1 - \sol{P}_2P_1f_1$. Similarly, $\sol{P}_1D_2f = f_2 - \sol{P}_1P_2f_2$. 

Consequently, we have 
\begin{align}\label{eq:RHS-P2D1-identity}
    &\frac{1}{|\Omega|} \sum_{c \in \cL_1 \cap \cL_2} (\sol{P}_2D_1f)(c) \cdot (\sol{P}_1D_2f)(c) \notag \\ 
    &= \frac{1}{|\Omega|} \sum_{c \in \cL_1 \cap \cL_2} f_1(c)f_2(c) - \frac{1}{|\Omega|} \sum_{c \in \cL_1 \cap \cL_2}f_1(c) \cdot (\sol{P}_1P_2f_2)(c) \notag \\ &\qquad -\frac{1}{|\Omega|}\sum_{c \in \cL_1 \cap \cL_2}(\sol{P}_2 P_1f_1)(c)\cdot f_2(c) + \frac{1}{|\Omega|} \sum_{c \in \cL_1 \cap \cL_2} (\sol{P}_2 P_1f_1)(c) \cdot (\sol{P}_1P_2f_2)(c) \notag \\
    &\underset{\text{Claim~\ref{claim:ip-sum}}}{=} \frac{1}{|\Omega|} \sum_{c \in \cL_1 \cap \cL_2} f_1(c) f_2(c) + \langle f_1, P_1P_2f_2 \rangle_{\mu} + \langle P_2P_1f_1, f_2 \rangle_{\mu} +\frac{1}{|\Omega|} \sum_{c \in \cL_1 \cap \cL_2} (\sol{P}_2P_1f_1)(c) \cdot (\sol{P}_1P_2f_2)(c) \notag \\
    &\underset{P_1, \, P_2 \text{ self adjoint}}{=} \frac{1}{|\Omega|} \sum_{c \in \cL_1 \cap \cL_2} f_1(c)f_2(c) + 2 \langle P_1f_1, P_2 f_2 \rangle_{\mu} + \frac{1}{|\Omega|} \sum_{c \in \cL_1 \cap \cL_2} (\sol{P}_2P_1f_1)(c) \cdot (\sol{P}_1P_2f_2)(c).
\end{align}
where in the second equality we use \eqref{eq:f1c0} and that $\sol{P}_1P_2f_2$ can be naturally extended to a function on $L^2(\mu)$ by looking at $P_1P_2f_2.$
Combining \eqref{eq:RHS-P2D1-identity} and \eqref{eq:LHS-P2D1-identity} proves \eqref{eq:P2D1-rearrange}. 

\begin{claim}\label{claim:P2f2-D2f2}
    For $i\in\{1,2 \}$, we have 
    \[ \|P_if_i\|_{\mu} \leq \frac{\|D_if_i\|_{\mu}}{\sqrt{(\ell_1-1)(\ell_2-1)-1}}. \]
\end{claim}
\begin{proof}
    We first relate $\|P_1f_1 \|_{\mu}$ to $\|f_1 \|_{\mu}$. To that end, because $\sum_{c \in \cL_1}f_1(c) = 0$, we have the explicit form 
\begin{equation}\label{eq:explicit-P1}
     P_1f_1(c)=
    \begin{cases}
        -\frac{f_1(c)}{\ell_1-1},&c\in\cL_1\cap\cL_2,\\
        0,&c\in\cL_2\setminus\cL_1,
    \end{cases}
\end{equation}
As such, we may write 
\begin{align*}
    \|P_1f_1\|_{\mu}^2 \underset{\eqref{eq:explicit-P1}}{=} \frac{1}{|\Omega|(\ell_1-1)}\sum_{c \in \cL_1 \cap \cL_2} f_1(c)^2 \leq \frac{\|f_1\|_{\mu}^2}{(\ell_1-1)(\ell_2-1)}.
\end{align*}
where for the inequality we use $\PP_{\mu}[X_1 = c] = \frac{\ell_2- 1}{|\Omega|}$ for $c \in \cL_1 \cap \cL_2$.

In turn this implies that 
\[(\ell_1 -1)(\ell_2-1)\|P_1f_1\|_{\mu}^2 \leq \|f_1 \|_{\mu}^2 \underset{\text{Pythagoras}}{=} \|D_1 f_1 \|_{\mu}^2 + \|P_1f_1 \|_{\mu}^2 \]
which rearranges to the desired.
\end{proof}

\begin{claim}\label{claim:ip-sum}
    Let $g_1,g_2\in L^2(\mu)$ such that $g_1$ is function of only $x_1$ and $g_2$ a function of only $x_2$ be such that $\sum_{c \in \cL_1} g_1(c,\cdot)=0$.
    Then 
    \[ \langle g_1, g_2 \rangle_\mu = - \frac{1}{|\Omega|}\sum_{c \in \cL_1 \cap \cL_2} g_1(c,\cdot) g_2(\cdot,c).\]
    The same conclusion would also hold if we instead only had $\sum_{c \in \cL_2} g_2(\cdot, c) = 0 $. Furthermore, the same conclusion holds even if $g_1\in \R^{\cL_1}$ (or $g_2\in \R^{\cL_2}$) simply by extending $g_1(x_1,x_2)=g_1(x_1)$ for any $(x_1,x_2)\in\Omega$.
\end{claim}
\begin{proof}
    We may write 
    \begin{align*}
        \langle g_1, g_2 \rangle_{\mu} &= \frac{1}{|\Omega|} \bigg[ \bigg( \sum_{x \in \cL_1} g_1(x, \cdot)\bigg) \bigg( \sum_{y \in \cL_2} g_2(\cdot, y) \bigg) -\sum_{c \in \cL_1 \cap \cL_2} g_1(c, \cdot )g_2(\cdot, c)\bigg] \\
        &= - \frac{1}{|\Omega|} \sum_{c \in \cL_1 \cap \cL_2}g_1(c, \cdot )g_2(\cdot, c). \qedhere
    \end{align*}
\end{proof}
This completes the proof of Lemma~\ref{lem:edge-frame}.
\end{proof}

\section{Local Inequalities on Subgraphs under a Disallowed Color}
In this section, we extend the machinery of the previous section to subgraphs $H\subset G$ where a partial proper coloring on $G\smallsetminus H$ is fixed such that a specific color $c$ is disallowed to be used in $H$. Later we will see that the machinery we develop here allows us to study local inequalities when $H$ has a vertex connected to every other vertex (which we will call a \emph{cone}). Specifically, we will apply these inequalities to the leaves of such a cone. 

First, we introduce some notation. We define the set of proper colorings avoiding $c$ by
\[
 \Omega^{-c}_H
 =\left\{\bx\in \prod_{i \in H}\cL^{-c}_i :
x_i\neq x_j
      \text{ for every }ij \in E(H)\right\}.
\]
We define
$\mu_H^{-c}$ to be the uniform distribution over all proper colorings in $(\Omega_H^{-c})$. We write $L^2(\mu^{-c})$ for the space of real-valued functions on
$\Omega^{-c}$. 
The operators $P_i^H$, $D_i^H$, and $\Pi_{i,c}^H$ are those from the previous section, acting on $L^2(\mu_H)$, and we write
$D=\sum_{i\in H}D_i^H$. As before, $\bx^{i\to a}$ agrees with $\bx$ except that
its $i$th coordinate is $a$.

\begin{fact}\label{fact:disallowed-color}
For every $h\in L^2(\mu_H)$, every color $c$, every
$x\in\Omega_H^{-c}$, and every $i\in H$, we have
\[
\mathbb E_{a \sim \mathcal L_i^{-c-\bx_{N_H(i)}}}[h(\bx^{i \to a})]
-\mathbb E_{a \sim \mathcal L_i^{-\bx_{N_H(i)}}}[h(\bx^{i \to a})]
=(\Pi_{i,c}^Hh)(\bx).
\]
\end{fact}

\begin{proof}
If $c\notin\mathcal L_i$, both sides of the first identity are zero.
Otherwise, we write
    \begin{align*}
       \E_{a\sim \cL_i^{-c - \bx_{N_H(i)}}} h(\bx^{i \to a}) - \E_{a \sim \cL_i^{-\bx_{N_H(i)}}} h(\bx^{i \to a}) 
        &= -\frac{h(\bx^{i \to c}) - \frac{1}{\ell_i^{-\bx_{N_H(i)}}} \sum_{a \in \cL_i^{-\bx_{N_H(i)}}} h(\bx^{i \to a})}{ \ell_i^{-\bx_{N_H(i)}} - 1} \\
        &= e_{i,c}(\bx) \cdot  \frac{h(\bx^{i \to c})-\frac1{\ell_i^{-\bx_{N_{H}(i)}}}\sum_{a\in\mathcal \cL_i^{-\bx_{N_{H}(i)}}}h(\bx^{i \to a})}{\sqrt{\ell_i^{-\bx_{N_{H}(i)}}-1}} \\
        &= e_{i,c}(\bx) \cdot 
  \frac1{\ell_i^{-\bx_{N_{H}(i)}}}\sum_{a\in\mathcal \cL_i^{-\bx_{N_{H}(i)}}}e_{i,c}(\bx^{i \to a})h(\bx^{i \to a}) \\
  &= \Pi_{i,c} h(\bx). \qedhere
    \end{align*}
\end{proof}

\begin{lemma}\label{lem:ABF}
For every $h\in L^2(\mu)$ and any $c \in \bigcup_{i \in V(H)} \cL_i$ where $\ell_i \geq s + \deg_H(i)$ for $i \in H$, we have 
    \[\E_{\mu} \Bigg[\mathbb{I} \{ c \not \in \bX \} \bigg( \sum_i \Pi_{i,c} h (\bX)\bigg)^2 \Bigg]  \leq (1-p_{-}^*)\bigg\| \sum_i \Pi_{i,c}h \bigg\|_{\mu}^2\leq \beta_{|H|}\bigg\| \sum_i \Pi_{i,c}h \bigg\|_{\mu}^2.\] 
    where $p_{-}^*:=(\frac{s-1}{s})^{|H|}.$
\end{lemma}
The main observation is to  rule out the event that $\sum_{i} \Pi_{i,c}h(\bx)$ has (almost) all of its mass on $\bx$'s such that $c \not \in \bX$. The idea is to decompose $\mu$ into a mixture of probability measures under which each $ \Pi_{i,c}h(\bx)$ has mean zero and the event $c\notin\bX$ consists of a single outcome
of probability at least $p_{-}^*$.

\begin{proof}
We rewrite $\mu$ in terms of resampling each vertex between its current
color and $c$, conditioned on the resulting coloring being proper.

We define
$\Omega^{-c}=\{\bx\in\Omega_H:c\notin\bx\}$.
We also write $p_{-c}=\P_\mu[c\notin\bX]>0$ and
$\mu^{-c}=\mu(\,\cdot\mid c\notin\bX)$. 

For $\bx\in\Omega^{-c}$, we define
\begin{equation}\label{eq:product-resampling}
 \nu_{c,\bx}(\by)
 =\prod_{i\in H}
 \left[
   \mathbb{I}\{y_i=x_i\}(1-p_{i,c}(\bx_{-i}))
   +\mathbb{I}\{y_i=c\}p_{i,c}(\bx_{-i})
 \right],
\end{equation}
and we define
$\pi_{c,\bx}
 =\nu_{c,\bx}(\,\cdot\mid\bY\in\Omega_H)$.
That is, we resample the coordinates independently according to
\eqref{eq:product-resampling} and condition on obtaining a proper coloring.
This conditioning has positive probability because $\bx\in\Omega_H$.
The event $c\notin\bY$ consists only of the outcome $\bY=\bx$.
Consequently, its probability satisfies
\begin{equation}\label{eq:no-change}
\begin{aligned}
 p_{-c}(\bx)
 :=\P_{\pi_{c,\bx}}(c\notin\bY)
   =\frac{\nu_{c,\bx}(\bx)}{\nu_{c,\bx}(\Omega_H)}
 \geq \prod_{i\in H}(1-p_{i,c}(\bx_{-i}))
  \geq \left(\frac{s-1}{s}\right)^{|H|}
  =p_{-}^*.
\end{aligned}
\end{equation}

We claim that, for every function $f:\Omega_H\to\mathbb R$,
\begin{equation}\label{eq:resampling-mixture}
 \E_\mu[f(\bY)]
 =p_{-c}\E_{\bX\sim\mu^{-c}}\left[
      \frac{1}{p_{-c}(\bX)}
      \E_{\bY\sim\pi_{c,\bX}}[f(\bY)]
   \right].
\end{equation}
For $\by\in\Omega_H$, define
$I_{\by}:=\{i\in H:y_i=c\}$.
Since $\by$ is a proper coloring, $I_{\by}$ is an independent set in $H$.
Next, we define $\nu_{\by}$ as the law of the random coloring $\bX$
obtained by independently sampling $X_i$ uniformly from
$\cL_i^{-\by_{N_H(i)}}\setminus\{c\}$
for each $i\in I_{\by}$, and setting $X_i=y_i$ for
$i\notin I_{\by}$.
For each possible outcome $\bx$ and every $i \in I_{\by}$, we have
$\bx_{N_H(i)}=\by_{N_H(i)}$ and so
\begin{align}\label{eq:reverse-resampling}
\nu_{\by}(\bx)
 &=\prod_{i\in I_{\by}}
   \frac{1}{\ell_i^{-\by_{N_H(i)}}-1}
   =\frac{\nu_{c,\bx}(\by)}{\nu_{c,\bx}(\bx)}
   =\frac{\pi_{c,\bx}(\by)}{p_{-c}(\bx)}.
\end{align}
For $\bx\in\Omega^{-c}$ outside the support of $\nu_{\by}$,
both $\nu_{\by}(\bx)$ and $\pi_{c,\bx}(\by)$ are zero.
Therefore, for $\by\in \Omega_H$ arbitrary, we have
\begin{align*}
 p_{-c}\E_{\bX\sim\mu^{-c}}\left[
    \frac{\pi_{c,\bX}(\by)}{p_{-c}(\bX)}
 \right] =\frac{1}{|\Omega_H|}
   \sum_{\bx\in\Omega^{-c}}\nu_{\by}(\bx) =\frac{1}{|\Omega_H|}=\mu(\by),
\end{align*}
where we used 
$p_{-c}=\frac{|\Omega^{-c}|}{|\Omega_H|}$. 
This proves \eqref{eq:resampling-mixture}.

Since $\bY=\bx$ when $c\notin\bY$ and writing $g_c(\bx) = \sum_i \Pi_{i,c} h(\bx)$, we have
\begin{align*}
 \E_{\bY\sim\mu}[\mathbb{I}\{c\notin\bY\}g_c(\bY)^2]
 &=p_{-c}\E_{\bX\sim\mu^{-c}}[g_c(\bX)^2]\\
 &=p_{-c}\E_{\bX\sim\mu^{-c}}\left[
   \frac{1}{p_{-c}(\bX)^2}
   \left(\E_{\pi_{c,\bX}}
     [\mathbb{I}\{c\notin\bY\}g_c(\bY)]\right)^2
   \right]\\
 &\underset{\forall i: \E_{\pi_{c,\bX}}[\Pi_{i,c}h|Y_{-i}]=0}{=}
 p_{-c}\E_{\bX\sim\mu^{-c}}\left[
   \frac{1}{p_{-c}(\bX)^2}
   \left(\E_{\pi_{c,\bX}}
     [(\mathbb{I}\{c\notin\bY\}-p_{-c}(\bX))g_c(\bY)]\right)^2
   \right]\\
 &\underset{\text{Cauchy--Schwarz}}{\leq}
p_{-c}\E_{\bX\sim\mu^{-c}}\left[
   \frac{1-p_{-c}(\bX)}{p_{-c}(\bX)}
   \E_{\pi_{c,\bX}}[g_c(\bY)^2]
   \right]\\
 &\underset{\eqref{eq:no-change}}{\leq}
 (1-p_{-}^*)p_{-c}\E_{\bX\sim\mu^{-c}}\left[
   \frac{1}{p_{-c}(\bX)}
   \E_{\pi_{c,\bX}}[g_c(\bY)^2]
   \right]\\
 &\underset{\eqref{eq:resampling-mixture}}{=}
 (1-p_-^*)\E_\mu[g_c(\bX)^2]. \qedhere
\end{align*}
\end{proof}

\section{Local Inequalities on Cones}

In this section we prove Theorem \ref{thm:local}. Before doing that we set up some notation.

Let $H= \{ v\} \cup H'$ be a cone centered at $v$.  We write
\[
    \Omega_{H'}:=\{\bx_{H'}:\bx\in\Omega\} \quad\text{and}\quad
    \Omega_{H'}^{-c}:=\{\bx\in\Omega_{H'}:c\notin\bx\}.
\]

We will often abuse notation and identify a leaf function $h:\Omega_{H'}\to\mathbb R$ 
with $(c,\bx)\mapsto h(\bx)$ on $\Omega$.
Once we pin $\bx_{-i}$, we will write $h_i(a)=h(\bx_{-i},a)$.

Let $\mu_{H'}^{-v}$
be the uniform distribution on the proper colorings of leaves from $\prod_i \cL_i$ when we assume that the center is assigned a dummy color that does not belong to any of the leaf lists (note that vertices outside of $H$ are already pinned and the lists of $u_i$'s and $v$ are chosen based on that pinning).
Note that a uniformly random proper coloring of $H$ can be obtained by first choosing $\bx$ with probability proportional to
$\ell_v^{-\bx}$ under $\mu_{H'}^{-v}$, and then choosing a uniform
color $c\in\cL_v^{-\bx}$.

We further define the operators $D_i^{-v,H'}$ and $\Pi_{i,c}^{-v,H'}$ acting on
$L^2(\mu_{H'}^{-v})$ to denote  operators obtained by
pinning the center to the dummy color.

Throughout the section, unless otherwise specified, we write $D_u = D_u^H$ and $P_u = P_u^H$ for $u \in H$. We also write,
$D_{H'}=\sum_{i\in H'}D_i$. 
For a cone $H$, we say a function $h:\Omega_H\to\R$ is a leaf function if it does not depend on the color  of the center $v$. 
Note that if $h$ is a leaf function, $D_{H'}h$ may depend on the center
color. 
\localCorollary*

\begin{proof}[Proof of Theorem~\ref{thm:local}]
In the first step we ``decouple'' the action of leaf update $D_{H'}$  on a leaf function versus a non-leaf function. To do that, 
    we use the following leaf function $h^*$.
\begin{equation}\label{eq:objective-new}
       h^* := \argmin_{h \text{ leaf function}}  \langle f-h, D_{H'}(f-h) \rangle_\mu
        = \argmin_{h \text{ leaf function}} \sum_{i \in H'} \|D_i(f-h)\|_\mu^2.
    \end{equation}
    Letting $g = f - h^*$, we will show towards the end of the proof that $D_{H'}g$ is orthogonal to all leaf functions. In particular, 
    \begin{equation}\label{eq:Dug-orthogonal-new}
        P_vD_{H'}g = 0.
    \end{equation}
    Having this, it follows that
    \[
        P_vD_{H'}h^* = P_vD_{H'}f \quad\text{and}\quad
        \langle h^*,D_{H'}g\rangle_\mu \underset{h^* \text{ leaf function}}{=} \langle P_vh^{*},D_{H'}g \rangle \underset{\eqref{eq:Dug-orthogonal-new}}{=}0,
    \]
    and consequently
    \begin{equation}\label{eq:f-h*}
        \langle f,D_{H'}f\rangle_\mu
        = \langle h^*,D_{H'}h^*\rangle_\mu
          + \langle g,D_{H'}g\rangle_\mu
        \geq \langle h^*,D_{H'}h^*\rangle_\mu.
    \end{equation}

Having decoupled $g,h^*$, we write,
    \begin{align}\label{eq:Dvf-DH'f}
        2\langle D_vf,D_{H'}f\rangle_\mu
        &\underset{\substack{D_v^2 = D_v \\
            D_v\text{ is self-adjoint}}}{=}
        2\langle D_vf,D_vD_{H'}h^*\rangle_\mu
        + 2\langle g,D_vD_{H'}g\rangle_\mu
        \notag \\
        &\underset{\eqref{eq:Dug-orthogonal-new}}{=}
        2\langle D_vf,D_vD_{H'}h^*\rangle_\mu
        + 2\langle g,D_{H'}g\rangle_\mu.
        \end{align}
 Since $D_{H'}$ is PSD (w.r.t. $\mu$) we can drop the second term. So, the main goal is to analyze $\langle D_vf, D_vD_{H'}h^*\rangle_\mu$. Let $\alpha:=\frac{\beta_d}{s-1}$. First, we use a careful AM-GM to de-couple $f$, 
        \begin{align*}
        2 (1-\alpha) \cdot  \langle D_vf,D_{H'}f\rangle_\mu
        &\geq
        2 \langle D_vf, (1-\alpha) D_vD_{H'}h^* \rangle_\mu \geq
        -(1-\alpha) \bigg(t\|D_vf\|_\mu^2
        +
        \frac{1}{t} \|D_vD_{H'}h^* \|_\mu^2 \bigg),
        \end{align*}
        where the second inequality uses the fact that, for any two vectors $\bx,\by$, $2\langle \bx,\by\rangle
            \geq -t\|\bx\|^2-\frac{\|\by\|^2}{t}$.    
        So, the main difficulty is to upper-bound $\norm{D_vD_{H'}h^*}_\mu^2.$
        The observation is that if $D_{H'}h^*$ is still a leaf function then $D_vD_{H'}h^*=0$. Otherwise, it must have some dependency on the center. But, as we will see in Proposition \ref{prop:cone} (below) it is not really too far from being a leaf function.
        \begin{align*}
        &\hspace{-4em}
        \underset{\text{Proposition~\ref{prop:cone}}}{\geq}
        -(1-\alpha) t \|D_v f \|_{\mu}^2 - \frac{\alpha}{t} \|P_vD_{H'} h^* \|_{\mu}^2 - O\bigg(  \frac{\Delta}{ts^3} \bigg) \langle h^*, D_{H'} h^* \rangle_{\mu}
        \\
        &\hspace{-3em}
        \underset{P_vD_{H'}h^*=P_vD_{H'}f}{\geq}
        -(1-\alpha) t\|D_vf\|_\mu^2
        -
        \frac{\alpha}{t}
        \|P_v D_{H'} f\|_\mu^2
        -
        O\left(\frac{\Delta}{ts^3}\right)
        \langle f,D_{H'}f\rangle_\mu.
    \end{align*}
On the other hand, by another AM-GM on $\langle D_vf, D_vD_{H'}f\rangle_\mu$ we get
    \[ 2\alpha \langle D_vf, D_vD_{H'} f \rangle \geq - \alpha t  \|D_vf\|_{\mu}^2 - \frac{\alpha}{t}\|D_v D_{H'} f\|_{\mu}^2\]
Summing up the above two inequalities proves the theorem. 

It remains to  prove \eqref{eq:Dug-orthogonal-new}.
    Let $\eps \in \R$ and $\ell$ be an arbitrary leaf function.
    By definition of $h^*$ being the minimizer of
    \eqref{eq:objective-new}, it follows that
    \[
        0
        =
        \left.
        \frac{d}{d\eps}
        \langle f-h^*-\eps\ell,
                D_{H'}(f-h^*-\eps\ell)\rangle_\mu
        \right|_{\eps=0}
        =
        -2\langle\ell,D_{H'}(f-h^*)\rangle_\mu
        =
        -2\langle\ell,D_{H'}g\rangle_\mu.
    \]
    This shows that $D_{H'}g$ is orthogonal to all leaf functions.
    Consequently, it follows that for any function $h$, we have
    \[
        \langle (I-D_v)h,D_{H'}g\rangle_\mu
        =
        \langle h,(I-D_v)D_{H'}g\rangle_\mu
        =0.
    \]
    Since the choice of $h$ was arbitrary, $(I-D_v)D_{H'}g=0$.
\end{proof}

The main statement we needed in the above proof is the following proposition which we will prove in the rest of this section. It says that for any leaf function $h$, $Dh$ is almost a leaf function. Note that if $Dh$ were indeed a leaf function then $D_vDh=0$. 

\begin{proposition}\label{prop:cone}
Let $H = \{v \} \cup H'$ be a cone, and suppose that the list of colors at every vertex $u \in H$ satisfies $\ell_u \geq s + \deg_H(u)$ and $s \geq 5$. For every leaf function $h$, we have
\[
\bigg( 1 - \frac{\beta_d}{s-1}\bigg) \|D_vD_{H'}h\|_\mu^2
\le\frac{\beta_d}{s-1}\|P_v D_{H'}h\|_\mu^2
+O\left(\frac{\Delta}{s^3}\right)
 \langle h,D_{H'}h\rangle_\mu,
\]
where $d = |H'|$.
\end{proposition}

We first construct a leaf function approximating $D_{H'}h$.

\begin{lemma}[$D_{H'}h$ is approximately a leaf function]
\label{lem:approx-leaf}
Let $H = \{v \} \cup H'$ be a cone and the list of colors at every vertex $u \in H$ is assumed to satisfy $\ell_u \geq s + \deg_H(u)$ and $s \geq 5$. For every leaf function $h$, there exists a leaf function $h'$ such that
\begin{equation}
\|h'-D_{H'}h\|_\mu^2
\le\frac{\beta_d}{s-1}\|h'\|_\mu^2
+O\left(\frac1{s^3}\right)
 \sum_{i\in H'}\deg_{H'}(i)\|D_i h\|_\mu^2,
\end{equation}
where $d = |H'|$.
\end{lemma}

Assuming Lemma~\ref{lem:approx-leaf}, we prove Proposition~\ref{prop:cone}. 

\begin{proof}[Proof of Proposition~\ref{prop:cone}]
We let $h'$ be the leaf function as in Lemma~\ref{lem:approx-leaf}. Since $D_vh'=0$ and
$D_v$ is an orthogonal projection, we have
\begin{align*}
    \norm{D_vD_{H'}h}_\mu^2
    &=
    \left\|
        D_v\left(
            D_{H'}h-
            \left(1-\frac{\beta_d}{s-1}\right)h'
        \right)
    \right\|_\mu^2\\
    &\leq
    \left\|
        D_{H'}h-
        \left(1-\frac{\beta_d}{s-1}\right)h'
    \right\|_\mu^2\\
    &=
    \frac{\beta_d}{s-1}\norm{D_{H'}h}_\mu^2+
    \left(1-\frac{\beta_d}{s-1}\right)
    \left[
        \norm{D_{H'}h-h'}_\mu^2
        -\frac{\beta_d }{s-1}\norm{h'}_\mu^2
    \right]\\
    &\underset{\text{Lemma~\ref{lem:approx-leaf}}}{\leq}
    \frac{\beta_d}{s-1}\norm{D_{H'}h}_\mu^2
    +O\left(\frac1{s^3}\right)
    \sum_{i\in H'}\deg_{H'}(i)\norm{D_i h}_\mu^2.
\end{align*}

Finally, each $D_i$ is an orthogonal projection and
$D_{H'}=\sum_{i\in H'}D_i$, so
\[
    \sum_{i\in H'}\deg_{H'}(i)\norm{D_i h}_\mu^2
    \leq
    \Delta\sum_{i\in H'}\norm{D_i h}_\mu^2
    =
    \Delta\langle h,D_{H'}h\rangle_\mu.
\]
Consequently,
\begin{align*}
    \norm{D_vD_{H'}h}_\mu^2
    &\leq
    \frac{\beta_d }{s-1}\norm{D_{H'}h}_\mu^2
    +O\left(\frac{\Delta}{s^3}\right)
      \langle h,D_{H'}h\rangle_\mu \\
      &\underset{\text{Pythagoras}}{=} \frac{\beta_d}{s-1} \Big(\norm{D_vD_{H'}h}_{\mu}^2 + \norm{P_vD_{H'}h}_{\mu}^2 \Big) + O\bigg( \frac{\Delta}{s^3} \bigg)\langle h, D_{H'}h\rangle_{\mu}.
\end{align*}
Rearranging proves the desired inequality.
\end{proof}

Next, before proving Lemma~\ref{lem:approx-leaf}, we compare
$\|D_i^{-v,H'}h\|_{\mu_{H'}^{-v}}^2$ and $\|D_i h\|_\mu^2$.

\begin{claim}\label{claim:leaf-cone}
If $\ell_v\geq s+|H'|$, then, for every leaf function $h$ and every $i\in H'$, we have
\[
\frac{|\Omega_{H'}|}{|\Omega|}\|D_i^{-v,H'}h\|_{\mu_{H'}^{-v}}^2
\le\frac1{s-1} \cdot \|D_i h\|_\mu^2.
\]
\end{claim}

\begin{proof}
Fix $i$ and a pinning $\bx_{-i}$ of the other leaves. We may write
\begin{align*}
\sum_{c\in\mathcal L_v^{-\bx_{-i}}}
 \sum_{a\in\mathcal L_i^{-c-\bx_{N_{H'}(i)}}}
 (D_i h)(c,\bx_{-i},a)^2
 &=\sum_{c \in \cL_v^{-\bx_{-i}}} \frac{1}{2 \ell_i^{-c - \bx_{N_{H'}(i)}}} \sum_{a,b \in \cL_i^{-c - \bx_{N_{H'}(i)}}}\Big( h\big(\bx^{i \to a} \big) - h\big(\bx^{i \to b} \big) \Big)^2 \\
&\geq
 \frac12\sum_{a,b\in\mathcal L_i^{-\bx_{N_{H'}(i)}}}
       (h(\bx^{i \to a})-h(\bx^{i \to b}))^2
 \sum_{c\in\mathcal L_v^{-a-b-\bx_{-i}}}
       \frac1{\ell_i^{-\bx_{N_{H'}(i)}}}\\
&\underset{\ell_v \geq s + \deg_H(v)}{\ge}
 \frac{s-1}{2\ell_i^{-\bx_{N_{H'}(i)}}}
 \sum_{a,b\in\mathcal L_i^{-\bx_{N_{H'}(i)}}}(h(\bx^{i \to a})-h(\bx^{i \to b}))^2\\
&=
 (s-1)\sum_{a\in\mathcal L_i^{-\bx_{N_{H'}(i)}}}
             (D_i^{-v,H'}h)(\bx_{-i},a)^2.
\end{align*}
In the expression on the LHS we use the vector $(c,\bx_{-i},a)$ to indicate that $v$ is colored $c$, vertex $i$ is colored $a$ and the rest are colored with $\bx_{-i}.$
Summing over $\bx_{-i}$ gives
\[
|\Omega|\|D_i h\|_\mu^2
\ge(s-1)|\Omega_{H'}|\|D_i^{-v,H'}h\|_{\mu_{H'}^{-v}}^2.
\]
Rearranging proves the desired claim.
\end{proof}
\begin{proof}[Proof of Lemma~\ref{lem:approx-leaf}]
We claim that $h':=\sum_{i\in H'}D_i^{-v,H'}h$ works.
For $(c,\bx)\in\Omega$, we have
\begin{equation}\label{eq:ell-DH}
    h'(\bx)-(D_{H'}h)(c,\bx) = \sum_{i \in H'} \big(\big(D_i^{-v,H'} - D_i)h\big)(c,\bx)
\underset{\text{Fact~\ref{fact:disallowed-color}}}{=}\begin{cases}\sum_{i\in H'}(\Pi_{i,c}^{-v,H'}h)(\bx) &c\notin\bx\\ 0&\text{otherwise.}
\end{cases}
\end{equation}
Note that in the last identity we use that if $c\notin\bx$, then $D_ih(c,\bx)$ acts as if the center is assigned a dummy color so it is same as $D_i^{-v,H'}.$

The lists on $H'$ satisfy
$\ell_i\ge s+1+\deg_{H'}(i)$, so we have 
\begin{align}\label{eq:intermediate-approx-leaf}
\|h'-D_{H'}h\|_\mu^2
&\underset{\eqref{eq:ell-DH}}{=}\frac1{|\Omega|}
  \sum_{c\in\mathcal L_v}\sum_{\bx\in\Omega_{H'}^{-c}}
  \left(\sum_{i\in H'}(\Pi_{i,c}^{-v,H'}h)(\bx)\right)^2 \notag \\
&= \frac{|\Omega_{H'}|}{|\Omega|} \sum_{c \in \cL_v} \E_{\mu_{H'}^{-v}}\bigg[ \mathbb{I}\{c \not \in \bX\} \bigg( \sum_{i \in H'} (\Pi_{i,c}^{-v,H'}h)(\bX) \bigg)^2 \bigg] \notag \\
&\underset{\text{Lemma~\ref{lem:ABF}}}{\le}\frac{\beta_d|\Omega_{H'}|}{|\Omega|}
  \sum_{c\in\mathcal L_v}
  \left\|\sum_{i\in H'}\Pi_{i,c}^{-v,H'}h\right\|_{\mu_{H'}^{-v}}^2 \notag \\
&\le\frac{\beta_d|\Omega_{H'}|}{|\Omega|}
  \sum_c\left\|\sum_{i\in H'}\Pi_{i,c}^{-v,H'}h\right\|_{\mu_{H'}^{-v}}^2 \notag \\
&\underset{\text{Lemma~\ref{lem:projection-sum}}}{\leq}  \frac{\beta_d |\Omega_{H'}|}{|\Omega|}  \bigg(\frac{s}{s-1} \| D_{H'}^{-v,H'}h\|_{\mu_{H'}^{-v}}^2+\frac{O(1)}{s^2}\sum_{i\in H'}\deg_{H'}(i)
       \|D_i^{-v,H'}h\|_{\mu_{H'}^{-v}}^2 \bigg) \notag \\
&= \frac{\beta_d |\Omega_{H'}|}{|\Omega|}  \bigg(\frac{s}{s-1} \|h'\|_{\mu_{H'}^{-v}}^2+\frac{O(1)}{s^2}\sum_{i\in H'}\deg_{H'}(i)
       \|D_i^{-v,H'}h\|_{\mu_{H'}^{-v}}^2 \bigg),
\end{align}
where in the last inequality we apply Lemma~\ref{lem:projection-sum} to $\Pi^{-v,H'}_{i,c}$ instead of $\Pi^{H'}_{i,c}$; but the lemma still goes through because we can assume  the center is pinned to any color including a dummy color. 

To simplify the above, we observe that every leaf coloring has at least $s$ extensions to the center, so that 
\begin{equation}\label{eq:Omega-ell}
    |\Omega|\|h'\|_\mu^2
=\sum_{\bx\in\Omega_{H'}}\ell_v^{-\bx}h'(\bx)^2
\underset{\ell_v\geq s + |H'|}{\ge} s|\Omega_{H'}|\|h'\|_{\mu_{H'}^{-v}}^2.
\end{equation}
Consequently, continuing from \eqref{eq:intermediate-approx-leaf} we obtain
\begin{align*}
\|h'-D_{H'}h\|_\mu^2
&\underset{\text{Claim~\ref{claim:leaf-cone}}, \, \eqref{eq:Omega-ell}}{\le}\frac{\beta_d}{s-1}\|h'\|_\mu^2
 +\frac{O(\beta_d)}{s^2(s-1)}
  \sum_{i\in H'}\deg_{H'}(i)\|D_i h\|_\mu^2.
\end{align*}
This proves Lemma~\ref{lem:approx-leaf}.
\end{proof}

\section{Analyzing Graphs ``far" From Bipartite}
In this section, we let $Q=\bigcup_{i\in G}\cL_i$.  We define 
\[ (\nabla_{i,c} f)(\bx) = \begin{cases}  f(\bx) - f(\bx^{i \to c}) &\text{if } c \in \cL_i^{-x_i - \bx_{N(i)}},  \\ 0 &\text{otherwise,} \end{cases} \]
and write $\nabla_i = (\nabla_{i,c})_{c \in Q}$. For every $f \in L^2(\mu)$ we have that $\nabla_{i,c}f \in L^2(\mu)$. In an abuse of notation we write 
\[\langle\nabla_i f,\nabla_j g\rangle_{Q \otimes \mu}
    := \sum_{c\in Q}
        \langle\nabla_{i,c}f,\nabla_{j,c}g\rangle_\mu.\] 

We define  ``unweighted'' Laplacian operators
\begin{equation}\label{eq:solD-i}
    \sol{D}_i f(\bx):=\ell_{i}^{-\bx_{N(i)}} D_if(\bx) =  \sum_{c\in \cL_i^{-\bx_{N(i)}}}
       \nabla_{i,c}f(\bx), \quad\text{and}\quad
    \sol{D}=\sum_{i\in G}\sol{D}_i.  
\end{equation}

We now show that, like $D_i,P_i$, the operators $\sol{D}_i$ are self-adjoint and PSD. Indeed, we write 
\begin{align*}
    \langle f,\sol{D}_i g\rangle_\mu
    &= \mathbb{E}_{\bx\sim\mu}\bigg[
        \sum_{c\in \cL_i^{-\bx_{N_G(i)}}}
        f(\bx)\bigl(g(\bx)-g(\bx^{i\to c})\bigr)
        \bigg]\\
    &= \frac{1}{2}\mathbb{E}_{\bx\sim\mu}\bigg[
        \sum_{c\in \cL_i^{-\bx_{N_G(i)}}}
        \big(f(\bx)-f(\bx^{i\to c})\big)
        \big(g(\bx)-g(\bx^{i\to c})\big)
        \bigg]\\
    &= \frac{1}{2}\langle\nabla_i f,\nabla_i g\rangle_{Q\otimes \mu}
     = \langle\sol{D}_i f,g\rangle_\mu.
\end{align*}
The same equality shows that $\langle f, \sol{D}_if \rangle_{\mu} = \frac{1}{2} \| \nabla_i f \|_{Q \otimes \mu}^2 \geq 0$. We also note that 
\begin{equation}\label{eq:Di2-Di}
    \|\sol{D}_if \|_{\mu}^2 = \E_{\mu}\Big[\Big(\ell_i^{-\bx_{N(i)}}\Big)^2 (D_if)^2\Big] \geq s \cdot \E_{\mu}\Big[\ell_i^{-\bx_{N(i)}} (D_if)^2\Big] = s \langle f, \sol{D}_if \rangle_{\mu}.
\end{equation}
where the inequality uses that $\ell^{-\bx_{{N}(i)}}_i\geq s.$
The following is the main lemma of our proof which, as we discussed in the intro, is a strengthening of Lemma \ref{lem:mixed-form-lower-bound}. In particular, analogous to the machinery of Matrix trickledown \cite[Proposition IV.4]{ALO21}, instead of lower-bounding $\langle \sol{D}_if,{D}_jf\rangle_{\mu}$ by $1/s$ fraction of the norms of $\norm{\sol{D}_if}_{\mu},\norm{\sol{D}_jf}_{\mu}$ we lower-bound by the inner-product of the gradient vectors.
\begin{lemma}\label{lem:pairwise}
    For a pair of vertices $i,j$, we have
    \[
        \langle\sol{D}_i f,\sol{D}_j f\rangle_{\mu}\geq
        \begin{cases}
            \langle\nabla_i f,\nabla_j f\rangle_{Q \otimes \mu},
                & \{i,j\}\in E,\\
            0, & \text{otherwise}.
        \end{cases}
    \]
\end{lemma}
This lemma immediately gives our main theorem of the section. 

\mainlambdamin*
\begin{proof}
Let $\lambda=|\lambda_{\min}(A_G)|$. We write
\begin{align*}
    \|\sol{D}f\|_{\mu}^2
    &\underset{\eqref{eq:Di2-Di}, \, \text{Lemma~\ref{lem:pairwise}}}{\geq} s \sum_i \langle f, \sol{D}_if \rangle_{\mu} +\sum_{i\neq j\in G}(A_G)_{ij}
            \langle\nabla_i f,\nabla_j f\rangle_{Q\otimes \mu}\\
    &\geq (s - 2\lambda) \sum_i \langle f, \sol{D}_if \rangle_{\mu}+ \lambda \sum_{i}\|\nabla_i f \|_{Q \otimes \mu}^2 + \sum_{i\neq j\in G}(A_G)_{ij}
            \langle\nabla_i f,\nabla_j f\rangle_{Q\otimes \mu} \\
    &=  (s-2\lambda)\sum_i \langle f, \sol{D}_if \rangle_{\mu}+\sum_{i,j\in G}\big(A_G+\lambda I\big)_{ij}        \langle \nabla_if, \nabla_jf \rangle_{Q       \otimes \mu}\geq (s-2\lambda)\sum_i \langle f, \sol{D}_if \rangle_{\mu} 
\end{align*}
where in the last inequality we use Lemma \ref{lemma:Hoeffding} and also that $A_G+\lambda I$ is a PSD matrix. This shows that the spectral gap of the Glauber dynamics is at least
$(s-2\lambda)/(n \cdot \max_{i} \ell_i)$, which implies the mixing time bound.
\end{proof}

It remains to prove Lemma \ref{lem:pairwise} which we do in the rest of this section.

\begin{proof}[Proof of Lemma~\ref{lem:pairwise}]
    If $\{i,j\}\notin E$, then $\sol{D}_i$ and $\sol{D}_j$ are
    commuting PSD operators, so
    $\langle\sol{D}_i f,\sol{D}_j f\rangle_{\mu}\geq0$.

Suppose now that $\{i,j\}\in E$. As in the proofs of the lemmas in Section~\ref{sec:pairlocalineq} we pin all vertices of $G$ but $i,j$ and we assume $\mu_{\{i,j\}}$ is the resulting uniform distribution over proper colorings $\Omega_{\{i,j\}}$ of $i,j$ and $\cL_i^{\{i,j\}},\cL_j^{\{i,j\}}$ are the resulting lists (we remove subscripts and superscripts for clarity).

Now, we note that 
\begin{align*}
    \langle \sol{D}_if, \sol{D}_jf\rangle_{\mu}         &\underset{\eqref{eq:solD-i}}{=} \Big\langle \sum_{c \in \cL_i} \nabla_{i,c}f, \sum_{c \in \cL_j} \nabla_{j,c}f \Big\rangle_{\mu} 
        = \langle \nabla_i f, \nabla_jf\rangle_{Q \otimes \mu} 
        +\frac{1}{4} \E_{\bx\sim\mu} \sum_{c_i,c_j} (\Box_{c_i,c_j} f(x_i, x_j))^2,
\end{align*}
where the inner-sum in the last expression is over all pairs of colors $c_i\in \cL_i^{-x_i-x_j}, c_j\in \cL_j^{-x_i-x_j}$ such that $c_i\neq c_j$ and
 $\Box_{c_i,c_j}f(x_i,x_j) = f(c_i,c_j) - f(c_i,x_j) - f(x_i, c_j) + f(x_i,x_j)$. 
 To see the last identity,consider a quadruple  $(x_i,x_j,c_i,c_j)$ occurring in the cross-terms and write
\[
x=f(x_i,x_j),\qquad y=f(c_i,x_j),\qquad z=f(x_i,c_j) \quad\text{and}\quad w=f(c_i,c_j).
\]
The four associated terms in the expanded sum are
\begin{align*}
(x-y)(x-z)+(y-x)(y-w) +(z-w)(z-x)+(w-z)(w-y) =(w-y-z+x)^2 =\big(\Box_{c_i,c_j}f(x_i,x_j)\big)^2. 
\end{align*}
Each such term is counted four times. \qedhere
\end{proof}

\bibliographystyle{alpha}
\bibliography{ref.bib}

\end{document}